\documentclass[twocolumn]{autartweb}    

\usepackage{graphicx}          
\usepackage{float}                               
\usepackage{cite}

\usepackage{amsmath}
\usepackage{amssymb}
\usepackage{gensymb}
\usepackage{mathtools}
\usepackage{upgreek} 
\usepackage{algorithm} 
\usepackage[noend]{algpseudocode} 
\usepackage{epstopdf} 
\usepackage{enumitem} 
\usepackage{color} 
\usepackage{hyperref}
\usepackage[normalem]{ulem} 

\newtheorem{problem}{Problem}
\newtheorem{remark}{Remark}
\newtheorem{assumption}{Assumption}
\newtheorem{definition}{Definition}
\newtheorem{lemma}{Lemma}
\newtheorem{corollary}{Corollary}

\DeclareMathOperator*{\argmin}{argmin}

\begin{document}

\begin{frontmatter}

\title{A Sequential Learning Algorithm for Probabilistically Robust Controller Tuning \thanksref{footnoteinfo}} 

\thanks[footnoteinfo]{Corresponding author R.~Chin.}

\author[unimelb,birmingham]{Robert Chin}\ead{r.chin4@student.unimelb.edu.au},    
\author[unimelb]{Chris Manzie},               
\author[anu]{Iman Shames},  
\author[unimelb]{Dragan Ne\v{s}i\'{c}},
\author[birmingham,ati]{Jonathan E. Rowe}

\address[unimelb]{Department of Electrical and Electronic Engineering, The University of Melbourne, Australia}  
\address[birmingham]{School of Computer Science, University of Birmingham, United Kingdom}             
\address[anu]{College of Engineering \& Computer Science, Australian National University, Australia}
\address[ati]{The Alan Turing Institute, United Kingdom}        
          
\begin{keyword}                           
Randomized algorithms; Robust control; Control of uncertain systems; Statistical learning theory; Ordinal optimization.               
\end{keyword}                             

\begin{abstract}                          
We introduce a sequential learning algorithm to address a robust controller tuning problem, which in effect, finds (with high probability) a candidate solution satisfying the internal performance constraint to a chance-constrained program which has black-box functions. The algorithm leverages ideas from the areas of randomised algorithms and ordinal optimisation, and also draws comparisons with the scenario approach; these have all been previously applied to finding approximate solutions for difficult design problems. By exploiting statistical correlations through black-box sampling, we formally prove that our algorithm yields a controller meeting the prescribed probabilistic performance specification. Additionally, we characterise the computational requirement of the algorithm with a probabilistic lower bound on the algorithm's stopping time. To validate our work, the algorithm is then demonstrated for tuning model predictive controllers on a diesel engine air-path across a fleet of vehicles. The algorithm successfully tuned a single controller to meet a desired tracking error performance, even in the presence of the plant uncertainty inherent across the fleet. Moreover, the algorithm was shown to exhibit a sample complexity comparable to the scenario approach.
\end{abstract}

\end{frontmatter}

\section{Introduction}

When there is probabilistic plant uncertainty, the use of randomised algorithms (RA) is a well-established technique for finding approximate solutions to otherwise difficult computational problems \cite{Tempo2013}. Some key results from this area of probabilistic robust control pertain to the sample complexity, i.e. the number of simulations needed for the RA to perform analysis or design until desired probabilistic specifications are met. In control analysis, the techniques used to obtain sample complexities were originally pioneered in \cite{Khargonekar1996, Tempo1996}, and based on well-known concentration inequalities such as the Chernoff bound. In control design, the first sample complexities for RA were provided in \cite{Vidyasagar1998}, using results from the paradigm of empirical risk minimisation in statistical learning theory (namely, using the Vapnik-Chervonenkis dimension). These randomised algorithms and variants thereof have seen numerous control applications. In fault detection, randomised algorithms have been used for finding a solution which ensures a low false alarm rate with high confidence \cite{Ding2019}. In \cite{Alpcan2003}, randomised algorithms were used for the estimation and analysis for the probability of stability in high speed communication networks.

In another line of literature, the premise of ordinal optimisation (OO) is to find approximate solutions to difficult stochastic optimisation problems \cite{Ho2007}. Introduced in \cite{Ho1992} for the optimisation of discrete-event dynamic systems, ordinal optimisation primarily operates under two principles: firstly that comparison by order (as opposed to comparison based on numerical difference) is more `robust' against noise, and secondly by goal softening, we can improve our chances at finding a successful solution. These advantages of ordinal comparison and goal softening have been theoretically demonstrated in \cite{Xie1997} and \cite{Lee1999} respectively. Noteworthy applications of OO include finding an approximate solution to the Witsenhausen problem (an unsolved problem in nonlinear stochastic optimal control) \cite{Deng1999}, and reducing the computational burden for rare event simulation of overflow probabilities in queuing systems \cite{Ho1995}.

We observe several commonalities between methods in OO and methods in RA for controller tuning. On the surface, they both seek to find approximate solutions to difficult design problems that are rendered intractable due to uncertainty. Moreover, they both employ a philosphy which can be roughly summarised as ``randomly sample many candidate solutions, simulate their performances, and pick the best observed one''. In OO, the practice of selecting the best is called the `horse race' rule, and its optimality was formally shown in \cite{Yang2002}. Additionally, \cite{Vidyasagar2001} discusses that although this strategy is what seems intuitively to be the best thing to do (and what has been done for decades), the usage of RA is justified through rigorous sample complexity estimates. A further similarity shared by OO and RA is the notion of goal softening, which can be used to control the degree of sub-optimality for the obtained solution.

This apparent connection between OO and RA had been recognised and briefly touched on in \cite{Khargonekar1996, Ishii2008}, but as of yet, has not been fully explored in the literature. In this paper, we address a controller tuning/design problem based on both OO and RA. That is, we present a randomised algorithm for solving a control design problem to meet a desired probabilistic performance specification, and characterise the sample complexity using recent work obtained in ordinal optimisation for copula models \cite{Chin2021}. Copulas are useful for modelling the dependence structure in multivariate distributions \cite{Joe2014}, and our results are best suited (i.e. the least conservative) when the underlying copula is not too unfavourably far from a Gaussian copula. To obtain our bounds, we also employ concentration inequalities akin to those often used in analysis of RA \cite{Dubhashi2009}.

The guarantee provided by our algorithm is that the tuned controller will meet a nominal performance threshold with high probability, in the presence of plant uncertainty. This result can be loosely compared to that of the \textit{scenario approach} RA for robust control design \cite{Calafiore2006}. The main distinction between the approaches are set of assumptions being operated under. Wherever the scenario approach assumes linearity/convexity of particular functions, we allow for the functions to be black-box (e.g. the result of a closed-loop simulation), but require there to be statistical correlation of some sense within the performance indicators.

Our algorithm also uses a stopping rule, so that the sample complexity is not known in advance, but rather is a random variable induced by the randomness over each run of the algorithm. As the decision of whether to stop is learned from the algorithm by drawing sequential samples, we refer to our algorithm as a `sequential learning' algorithm. Another sequential learning algorithm also appeared in \cite{Koltchinskii2000}, which built upon the work of \cite{Vidyasagar1998} with less conservative sample complexities. Their algorithm is based on the Rademacher bootstrap technique. Stopping rules in RA were also studied in \cite{Fujisaki2007} for designing linear quadratic regulators, while \cite{Alamo2015} investigated another class of sequential algorithms. A stopping rule is also considered by \cite{Bayraksan2012} for solving stochastic programs, in which the algorithm stops when the computed confidence widths of estimated quantities become sufficiently small; this is similar to the nature of our algorithm.

This paper is organised as follows. In Section \ref{sec:prelim}, we state the problem formulation and introduce the OO success probability. In Section \ref{sec:lcb}, bounds are developed leading up to a lower confidence bound for the OO success probability. This is followed in Section \ref{sec:algorithm} by our sequential learning algorithm, which applies the lower confidence bounds from the preceding section. Additionally, we provide a lower bound for the distribution function of the stopping time of the algorithm. The algorithm is then specialised to probabilistically robust controller tuning in Algorithm \ref{alg:tuning}, for which we state and prove our main result in Theorem \ref{thm:tuning}. Lastly in Section \ref{sec:numerical}, we apply our algorithm to a numerical example, which considers probabilistically robust tuning of model predictive controllers on a diesel engine air-path for a fleet of vehicles.

\section{Preliminaries}
\label{sec:prelim}

\subsection{Notation}

Throughout this paper, $\mathbb{R}$ denotes the set of real numbers and $\mathbb{N}$ denotes the set of natural numbers. We let $\left\lfloor \cdot \right\rfloor$ and $\left\lceil \cdot \right\rceil$ denote the integer floor and ceiling operators respectively. If $Q$ is a matrix, then $Q \succ 0$ means that $Q$ is positive definite. The function $\exp\left(\cdot\right)$ is the exponential function, the logarithm $\log\left(\cdot\right)$ is taken to be the natural logarithm, while the cotangent function is denoted by $\cot\left(\cdot\right)$. The probability of an event and expectation operator are denoted by $\operatorname{Pr}\left(\cdot\right)$ and $\mathbb{E}\left[\cdot\right]$ respectively, with context as to which probability space provided in subscripts (whenever further context is required). Equality in law between random elements is denoted with the binary relation $\underset{\mathrm{st}}{=}$. The indicator variable for the event $\mathcal{A}$ is given by $\mathbb{I}_{\mathcal{A}}$. The standard Gaussian cumulative distribution function (CDF) is denoted by $\Phi\left(\cdot\right)$, while its inverse (i.e. quantile function) is denoted by $\Phi^{-1}\left(\cdot\right)$. The exponential distribution with rate parameter $\lambda$ is represented by $\operatorname{Exp}\left(\lambda\right)$. The abbreviation i.i.d. stands for mutually independent and identically distributed.

\subsection{Problem Setup}

Consider a measurable \textit{system performance function}
\begin{equation}
J\left(\psi,\theta\right):\Psi\times\Theta\to\mathbb{R},
\end{equation}
with controller parameter $\theta \in \Theta$ from a topological space $\Theta$, and uncertain plant parameter $\psi \in \Psi$ from a topological space $\Psi$. The uncertainty over $\psi$ is represented by some probability distribution $\mathcal{P}_{\psi}$ over $\Psi$. Without loss of generality, we use the convention that a lower $J$ indicates better performance. Also introduce the \textit{ordinal comparison function}
\begin{equation}
\overline{J}\left(\theta\right): \Theta \to \mathbb{R},
\end{equation}
which allows for any two controller parameters to be compared, without respect to the specific value for the plant parameter $\psi$. As a concrete example for $\overline{J}\left(\theta\right)$, we could take for instance $\overline{J}\left(\theta\right) = J\left(\overline{\psi}, \theta\right)$ for some nominal value $\overline{\psi}$, such as $\overline{\psi} = \mathbb{E}_{\mathcal{P}_{\psi}}\left[\psi\right]$. An alternative example is to average out the uncertainty by taking $\overline{J}\left(\theta\right) = \mathbb{E}_{\mathcal{P}_{\psi}}\left[J\left(\psi, \theta\right)\right]$, supposing this expectation can be evaluated.

Our aim is to find a controller $\theta^{*}$ so that the system will perform `well' with high probability in the presence of plant uncertainty. To this end, let $J^{*} \in \mathbb{R}$ denote a nominal performance threshold, which is used to benchmark the performance $J\left(\psi, \theta^{*}\right)$. We then introduce the following problem statement.
\begin{problem}
Given $\gamma \in \Gamma \subset \left[0, 1\right]$ and nominal performance threshold $J^{*} \in \mathbb{R}$, find a controller $\theta^{*}$ (possibly at random) such that
\begin{equation}
\operatorname{Pr}_{\psi, \theta^{*}}\left(J\left(\psi,\theta^{*}\right)\leq J^{*}\right) \geq 1 - \gamma.
\label{eq:problem_probability}
\end{equation}
\label{prob:tuning}
\end{problem}

\subsection{Relation to Chance-Constrained Programming}

A relation can be drawn between Problem \ref{prob:tuning} and existing RAs for approximately solving chance-constrained programs. Consider the chance-constrained program
\begin{equation}
\begin{aligned}\min_{\theta \in \Theta} & \quad \overline{J}\left(\theta\right) \\
\text{subject to } & \operatorname{Pr}_{\psi}\left(J\left(\psi,\theta\right)\leq J^{*}\middle|\theta\right)\geq1-\varepsilon,\\
\end{aligned}
\label{eq:chance_constrained_program}
\end{equation}
for given $\varepsilon \in \left(0, 1\right]$, with feasible set denoted
\begin{equation}
\Theta^{\star} = \left\{\theta \in \Theta: \operatorname{Pr}_{\psi}\left(J\left(\psi,\theta\right)\leq J^{*}\middle|\theta\right)\geq1-\varepsilon \right\}.
\end{equation}
Sadly, chance-constrained programs are usually intractable to solve exactly. The scenario approach \cite{Calafiore2006} is an RA which solves exactly an approximate version of \eqref{eq:chance_constrained_program}, using a finite number of samples for the plant uncertainty (the `scenarios'). Existing theoretical results provide the sample complexity for the number of scenarios such that the randomised candidate solution $\theta^{*}$ is feasible to the original problem with high probability, i.e.
\begin{equation}
\operatorname{Pr}\left(\theta^{*}\in\Theta^{\star}\right)\geq1-\eta
\label{eq:feasible_set_prob}
\end{equation}
for given $\eta \in \left(0, 1\right]$. This is referred to as a \textit{two-level of probability} statement, since it decouples the probability due to the RA and probability due to plant uncertainty. However, a \textit{one-level of probability} interpretation encompassing randomness over both the RA and plant uncertainty follows from
\begin{align}
\operatorname{Pr}_{\psi,\theta^{*}}\left(J\left(\psi,\theta^{*}\right)\leq J^{*}\right) &\geq \begin{aligned}[t] & \operatorname{Pr}_{\psi}\left(J\left(\psi,\theta^{*}\right)\leq J^{*}\middle|\theta^{*}\in\Theta^{\star}\right) \\
& \times \operatorname{Pr}\left(\theta^{*}\in\Theta^{\star}\right)
\end{aligned} \\
&\geq \left(1-\varepsilon\right)\left(1-\eta\right) \\
&> 1-\varepsilon-\eta.
\end{align}
Hence the one-level of probability specification in Problem \ref{prob:tuning} is addressed with $\gamma = \varepsilon + \eta$, and the interpretation is that the algorithm finds with high probability a candidiate solution satisfying the internal constraint. Furthermore, the sample complexity results for the scenario approach are derived under the assumptions that $\overline{J}\left(\theta\right)$ is linear in $\theta$ and $J\left(\psi,\theta\right)$ is convex in $\theta$ for any value of $\psi$. This ensures that the resulting scenario program is convex. Recent extensions to the scenario approach \cite{Grammatico2015, Esfahani2015} allow for degrees of non-convexity in the scenario program, but still maintain some of the core assumptions in $\overline{J}\left(\theta\right)$ and $J\left(\psi,\theta\right)$. A related approach known as the \textit{sample approximation} approach \cite{Luedtke2008} allows for non-linearity of $\overline{J}\left(\theta\right)$ and non-convexity of $J\left(\psi,\theta\right)$, but the corresponding sample complexity results are valid either when $\Theta$ is a finite set, or in the case when the $\psi$ term can be separated out from $J\left(\psi,\theta\right)$.

Our proposed algorithm for Problem \ref{prob:tuning} may also be used to find a candidate solution $\theta^{*}$ satisfying the internal constraint $J\left(\psi,\theta^{*}\right)\leq J^{*}$ to \eqref{eq:chance_constrained_program} with high probability. In this view, our algorithm imposes less restrictive structure than existing algorithms, because our theoretical results can apply when $J\left(\psi,\theta\right)$ and $\overline{J}\left(\theta\right)$ are black-box functions (e.g. the result of a closed-loop simulation). However, we work with a qualitatively different set of assumptions: exploiting when $J\left(\psi,\theta\right)$ and $\overline{J}\left(\theta\right)$ share some positive statistical correlation (e.g. when they are related performance indicators). Intuitively, if a candidate solution $\theta$ is found such that $\overline{J}\left(\theta\right)$ is low, this is correlated with low $J\left(\psi,\theta\right)$, meaning that $\theta$ is more likely to satisfy the performance constraint. We use copulas to express the notion of correlation/dependence, discussed over the following subsection.

\subsection{Copula Modelling}

Several algorithms in OO and RA require a mechanism $\mathcal{P}_{\theta}$, which to randomly sample a candidate controller $\theta_{i} \in \Theta$. In order to find our candidate solution $\theta^{*}$ to Problem \ref{prob:tuning}, we propose applying the ``randomly sample and select the best to test" methodology by first sampling $n$ i.i.d. $\theta_{i} \sim \mathcal{P}_{\theta}$ and letting
\begin{equation}
\theta^{*} = \argmin_{\theta_{i} \in \left\{\theta_{1}, \dots, \theta_{n}\right\}}\overline{J}\left(\theta_{i}\right).
\end{equation}
Note that since $\mathcal{P}_{\theta}$ is treated as arbitrary, it is not necessarily required to sample uniformly from $\Theta$. Instead, the practitioner may elect to use a mechanism which induces lower values of $\overline{J}\left(\theta\right)$. For example, each $\theta_{i}$ could be the result of an i.i.d. run of a randomised optimisation algorithm which is bespoke to the properties of $\Theta$.

Once $\theta^{*}$ has been obtained, a single `test' of the system yields the performance $J\left(\psi^{\sharp}, \theta^{*}\right)$, from an independently realised plant $\psi^{\sharp} \sim \mathcal{P}_{\psi}$. This test performance naturally predicates on how well the two random variables $\overline{J}\left(\theta_{i}\right)$ and $J\left(\psi^{\sharp}, \theta_{i}\right)$ are correlated, via their dependence on $\theta_{i}$.  A strong correlation should suggest that well-performing $\overline{J}\left(\theta_{i}\right)$ is highly indicative of well-performing $J\left(\psi^{\sharp}, \theta_{i}\right)$, thus we would reasonably anticipate the test $J\left(\psi^{\sharp}, \theta^{*}\right)$ to also perform well.

To formalise the concept of dependence between $\overline{J}\left(\theta_{i}\right)$ and $J\left(\psi^{\sharp}, \theta_{i}\right)$, we use copulas \cite{Joe2014}. A copula is a multivariate distribution with uniform $\left(0, 1\right)$ marginals, so that via the inverse probability integral transform, a multivariate distribution may be represented with just its marginal distributions and a copula. A common choice for a copula model is the Gaussian copula, defined in the bivariate case as follows.
\begin{definition}[Bivariate Gaussian copula]
Let $\left(Z, X\right)$ be a bivariate standard Gaussian with correlation $\rho \in \left[-1, 1\right]$, i.e.
\begin{equation}
\begin{bmatrix}Z\\
X
\end{bmatrix}\sim\mathcal{N}\left(\begin{bmatrix}0\\
0
\end{bmatrix},\begin{bmatrix}1 & \rho\\
\rho & 1
\end{bmatrix}\right).
\end{equation}
Then the Gaussian copula with correlation $\rho$ is the distribution of $\left(\Phi\left(Z\right), \Phi\left(X\right)\right)$.
\label{def:bivariate_gaussian_copula}
\end{definition}
For a Gaussian copula, the dependence is neatly summarised with the correlation parameter $\rho$. However, not every family of copula will be parametrised with a correlation. Instead, a well-defined notion of correlation valid for any bivariate distribution is the Kendall correlation.
\begin{definition}[Population Kendall correlation]
For a bivariate distribution $\left(Z, X\right)$, the population Kendall correlation is defined as
\begin{equation}
\kappa = \mathbb{E}\left[\operatorname{sign}\left(Z - \grave{Z}\right)\operatorname{sign}\left(X - \grave{X}\right)\right],
\end{equation}
where $\left(\grave{Z}, \grave{X}\right)$ is an independent copy of $\left(Z, X\right)$.
\end{definition}
In this paper, it will be convenient to associate every bivariate distribution with a bivariate Gaussian copula, which we do so through the Kendall correlation.
\begin{definition}[Associated Gaussian copula]
For any bivariate distribution $\left(Z, X\right)$ with population Kendall correlation $\kappa$, the Gaussian copula associated with this distribution is defined as the bivariate Gaussian copula with correlation $\rho = \sin\left(\pi\kappa/2\right)$.
\label{def:associated_gaussian_copula}
\end{definition}
The formula $\rho = \sin\left(\pi\kappa/2\right)$ is from the relation between $\kappa$ and $\rho$ for a Gaussian copula \cite[Equation (9.11)]{Kendall1990}. As such, any bivariate distribution with a Gaussian copula has its own copula as the associated Gaussian copula.

\subsection{Standing Assumptions}

We are ready to list the standing assumptions of the paper, for which the main results are based on.
\begin{assumption}
The bivariate distribution for the performances $\left(\overline{J}\left(\theta_{i}\right), J\left(\psi^{\sharp}, \theta_{i}\right)\right)$ is continuous, and has population Kendall correlation $\kappa > 0$, however the value of $\kappa$ itself is unknown.
\label{assump:copula}
\end{assumption}
\begin{remark}
By Sklar's theorem \cite[Theorem 1.1]{Joe2014}, the continuity property in Assumption \ref{assump:copula} ensures that $\left(\overline{J}\left(\theta_{i}\right), J\left(\psi^{\sharp}, \theta_{i}\right)\right)$ has a unique copula.
\end{remark}
We also assume the following bound between the copula of $\left(\overline{J}\left(\theta_{i}\right), J\left(\psi^{\sharp}, \theta_{i}\right)\right)$, and its associated Gaussian copula.
\begin{assumption}
Let $\left(\widetilde{Z}, \widetilde{X}\right)$ denote the copula of $\left(\overline{J}\left(\theta_{i}\right), J\left(\psi^{\sharp}, \theta_{i}\right)\right)$ and let $\left(\widetilde{Z}, \widetilde{X}'\right)$ denote the Gaussian copula associated with $\left(\overline{J}\left(\theta_{i}\right), J\left(\psi^{\sharp}, \theta_{i}\right)\right)$. For a given $\nu \in \left[0, 1\right)$, then for all $z \in \left(0, 1\right)$ we have
\begin{equation}
\sup_{x \in \left(0, 1\right)}\left\{\operatorname{Pr}\left(\widetilde{X}' \leq x \middle|\widetilde{Z} = z\right) - \operatorname{Pr}\left(\widetilde{X} \leq x \middle|\widetilde{Z} = z\right)\right\} \leq \nu.
\label{eq:conditional_KS1}
\end{equation}
\label{assump:conditional_KS1}
\end{assumption}
\begin{remark}
Although our results apply to any distribution minimally satisfying Assumption \ref{assump:copula}, the condition \eqref{eq:conditional_KS1} in Assumption \ref{assump:conditional_KS1} is intuitively saying that the copula of $\left(\overline{J}\left(\theta_{i}\right), J\left(\psi^{\sharp}, \theta_{i}\right)\right)$ is not too unfavourably `far' (i.e. upper bounded by $\nu$) from that of a Gaussian copula. To elaborate further, a given bound on the Kolmogorov-Smirnov distance (i.e. supremum norm) or total variation distance \cite[\S 5.9]{Gray2011} between $\operatorname{Pr}\left(\widetilde{X}' \leq x \middle|\widetilde{Z} = z\right)$ and $\operatorname{Pr}\left(\widetilde{X} \leq x \middle|\widetilde{Z} = z\right)$ will imply \eqref{eq:conditional_KS1}. Moreover, if $\left(\overline{J}\left(\theta_{i}\right), J\left(\psi^{\sharp}, \theta_{i}\right)\right)$ is assumed to have a Gaussian copula, then \eqref{eq:conditional_KS1} is satisfied with $\nu = 0$.
\label{rem:conditional_KS1}
\end{remark}
We also require the nominal performance threshold $J^{*}$ to be feasible, in the following sense.
\begin{assumption}
The nominal performance threshold $J^{*}$ satisfies
\begin{equation}
\operatorname{Pr}_{\psi^{\sharp}, \theta_{i}}\left(J\left(\psi^{\sharp}, \theta_{i}\right) \leq J^{*}\right) > 0.
\end{equation}
\label{assump:feasible}
\end{assumption}
Lastly, we can forego exact knowledge about the distributions of $\mathcal{P}_{\psi}$, $\mathcal{P}_{\theta}$, but the standing assumption is that they can at the very least be sampled from (e.g. via a computer simulation).
\begin{assumption}
Samples can be drawn i.i.d. from the distributions $\mathcal{P}_{\psi}$ and $\mathcal{P}_{\theta}$.
\label{assump:sampling}
\end{assumption}
As a consequence, we can produce an i.i.d. sample from the distribution of $\left(\overline{J}\left(\theta_{i}\right), J\left(\psi^{\sharp}, \theta_{i}\right)\right)$, which we denote
\begin{equation}
\left(\overline{J}\left(\theta_{1}\right), J\left(\psi_{1}, \theta_{1}\right)\right), \dots, \left(\overline{J}\left(\theta_{n}\right), J\left(\psi_{n}, \theta_{n}\right)\right).
\label{eq:bivariate_sample_performances}
\end{equation}

Based on these assumptions, we propose Algorithm \ref{alg:tuning} (in Section \ref{sec:algorithm}) to address Problem \ref{prob:tuning}, with the formal statement contained in Theorem \ref{thm:tuning}.

\subsection{Ordinal Optimisation}

We briefly review some results from \cite{Chin2021}, which studied a particular success probability related to ordinal optimisation.

Consider $n$ i.i.d. copies of $\left(Z_{i}, X_{i}\right)$ drawn from the distribution of $\left(Z, X\right)$, which is continuous and has population Kendall correlation $\kappa > 0$. We observe $Z_{1}, \dots, Z_{n}$, and order these observations from best to worst, denoted by $Z_{1:n} \leq \dots \leq Z_{n:n}$. The best $m$ are selected, given by $Z_{1:n} \leq \dots \leq Z_{m:n}$, with respective $X$-values denoted as $X_{\left\langle 1\right\rangle}, \dots, X_{\left\langle m\right\rangle}$, which are initially unobserved. More explicitly, we have selected the pairs $\left(Z_{1:n}, X_{\left\langle 1\right\rangle}\right), \dots, \left(Z_{m:n}, X_{\left\langle m\right\rangle}\right)$.
\begin{definition}[OO success probability]
The success probability is defined as
\begin{align}
p_{\mathrm{success}}\left(n, m, \alpha\right) &:= \operatorname{Pr}\left(\bigcup_{i = 1}^{m}\left\{X_{\left\langle i\right\rangle} \leq x_{\alpha}^{*}\right\}\right) \label{eq:ord-opt-defn} \\
&= \operatorname{Pr}\left( \min_{i \in \left\{1, \dots, m\right\}}\left\{X_{\left\langle i\right\rangle}\right\} \leq x_{\alpha}^{*}\right),
\end{align}
where $x_{\alpha}^{*}$ with $\alpha \in \left(0, 1\right]$ is the $100\alpha$ percentile of the distribution of $X$, i.e. $\operatorname{Pr}\left(X \leq x_{\alpha}^{*}\right) = \alpha$.
\label{def:ord_opt}
\end{definition}

Several properties in \cite{Chin2021} pertaining to the OO success probability can be specialised to the case of Gaussian copulas.

\begin{thm}[Gaussian OO success probability]
If the distribution of $\left(Z, X\right)$ in Definition \ref{def:ord_opt} has a Gaussian copula with correlation $\rho > 0$, then the ordinal optimisation success probability \eqref{eq:ord-opt-defn}, here denoted $p_{\mathrm{success}}^{\mathcal{N}}\left(n, m, \alpha, \rho\right)$, satisfies the following properties.
\begin{enumerate}[label=(\alph*)]
\item (Monotonicity in $m$) Given the triplet $\left(\bar{n}, \bar{\alpha}, \bar{\rho}\right) \in \mathbb{N} \times \left(0, 1\right) \times \left(0, 1\right)$, then
\begin{equation}
p_{\mathrm{success}}^{\mathcal{N}}\left(\bar{n}, m, \bar{\alpha}, \bar{\rho}\right) \leq p_{\mathrm{success}}^{\mathcal{N}}\left(\bar{n}, m', \bar{\alpha}, \bar{\rho}\right)
\end{equation}
for all $m \in \left[1, n\right)$ and $m' \in \left[m, n\right]$. \label{thm:OO_monotonicity_m}
\item (Monotonicity in $n$) Given the triplet $\left(\bar{m}, \bar{\alpha}, \bar{\rho}\right) \in \mathbb{N} \times \left(0, 1\right] \times \left(0, 1\right]$, then for all $n \leq n'$ such that $n' \geq \bar{m}$ and $n \in \left[\bar{m}, n'\right]$, we have
\begin{equation}
p_{\mathrm{success}}^{\mathcal{N}}\left(n, \bar{m}, \bar{\alpha}, \bar{\rho}\right) \leq p_{\mathrm{success}}^{\mathcal{N}}\left(n', \bar{m}, \bar{\alpha}, \bar{\rho}\right).
\end{equation}\label{thm:OO_monotonicity_n}
\item (Convergence of success probability) Given the triplet $\left(\bar{m}, \bar{\alpha}, \bar{\rho}\right) \in \mathbb{N} \times \left(0, 1\right] \times \left(0, 1\right]$, then
\begin{equation}
\lim_{n \to \infty}p_{\mathrm{success}}^{\mathcal{N}}\left(n, \bar{m}, \bar{\alpha}, \bar{\rho}\right) = 1.
\end{equation}\label{thm:approx-convergence-part-b}
\item (High probability of success) Given the triplet $\left(\bar{m}, \bar{\alpha}, \bar{\rho}\right) \in \mathbb{N} \times \left(0, 1\right] \times \left(0, 1\right]$, and for any $\delta \in \left(0, 1\right]$, there exists an $\widetilde{n}_{\delta}\left(\bar{\alpha}, \bar{\rho}\right) < \infty$ such that
\begin{equation}
p_{\mathrm{success}}^{\mathcal{N}}\left(n, \bar{m}, \bar{\alpha}, \bar{\rho}\right) \geq 1 - \delta,
\end{equation}
for all $n \geq \widetilde{n}_{\delta}\left(\bar{\alpha}, \bar{\rho}\right)$. \label{thm:high_probability_success}
\item (Lower bound for success probability) For any $\omega \in \left(0, \pi/2\right)$, let
\begin{gather}
c_{1} = \dfrac{1}{2} - \dfrac{\omega}{\pi} \\
c_{2} = \dfrac{\cot{\omega}}{\pi - 2\omega}
\end{gather}
and
\begin{gather}
\mu_{n} = -\sqrt{\dfrac{\log\left(nc_{1}\right)}{c_{2}}} \label{eq:mu-constructed} \\
\sigma_{n}^{2} = \dfrac{-\log\log 2}{2c_{2}\left(\log\left(nc_{1}\right) - \log\log 2\right)}. \label{eq:sigma-constructed}
\end{gather}
Then there exists some $n^{*}\left(\omega\right) \in \mathbb{N}$ such that for all $n \geq n^{*}\left(\omega\right)$, $m \in \left[1, n\right]$, $\rho \in \left(0, 1\right]$, $\alpha \in \left(0, 1\right]$, we have
\begin{align}
p_{\mathrm{success}}^{\mathcal{N}}\left(n, m, \alpha, \rho\right) &\geq \Phi\left(\dfrac{\Phi^{-1}\left(\alpha\right) - \rho\mu_{n}}{\sqrt{1 - \rho^{2} + \rho^{2}\sigma_{n}^{2}}}\right) \nonumber \\
&=: \widehat{p}_{\mathrm{success}, \omega}^{\mathcal{N}}\left(n, m, \alpha, \rho\right).
\label{eq:success-prob-copula-lb}
\end{align}
\label{thm:OO_success_prob_copula_lb}
\end{enumerate}
\label{thm:OO_properties}
\end{thm}
From Theorem \ref{thm:OO_properties}\ref{thm:OO_success_prob_copula_lb}, an optimised lower bound (optimised with respect to $\omega$) can be derived as
\begin{align}
p_{\mathrm{success}}^{\mathcal{N}}\left(n, m, \alpha, \rho\right) &\geq \sup_{\omega \in \Omega_{n}}\widehat{p}_{\mathrm{success}, \omega}^{\mathcal{N}}\left(n, m, \alpha, \rho\right) \nonumber \\
&=: \widehat{p}_{\mathrm{success}}^{\mathcal{N}}\left(n, m, \alpha, \rho\right),
\label{eq:success-prob-copula-lb-optimised}
\end{align}
where $\Omega_{n} \subset \left(0, \pi/2\right)$ is the set of all $\omega$ such that $n \geq n^{*}\left(\omega\right)$. A time complexity $O\left(1\right)$ numerical implementation of \eqref{eq:success-prob-copula-lb-optimised} is detailed in \cite{Chin2021}.

\section{Success Probability Lower Confidence Bound}
\label{sec:lcb}

Problem \ref{prob:tuning} can be framed in the context of OO, by taking
\begin{equation}
\left(Z_{i}, X_{i}\right) \underset{\mathrm{st}}{=} \left(\overline{J}\left(\theta_{i}\right), J\left(\psi^{\sharp}, \theta_{i}\right)\right)
\label{eq:substitution}
\end{equation}
in Definition \ref{def:ord_opt} with $m = 1$. However, a value for $\alpha = \operatorname{Pr}_{\psi^{\sharp}, \theta_{i}}\left(J\left(\psi^{\sharp}, \theta_{i}\right) \leq J^{*}\right)$ is not explicitly mentioned in Problem \ref{prob:tuning}, nor can the results in Theorem \ref{thm:OO_properties} be readily applied since $\left(Z, X\right)$ may generally not have a Gaussian copula. In this section, we overcome these obstacles by developing a lower confidence bound for the OO success probability. This is to be derived from lower confidence bounds for $\alpha$, and for $\rho$, the latter being the correlation of the associated Gaussian copula.

To facilitate this, we will work more abstractly with a continuous bivariate distribution for $\left(Z, X\right)$ with Kendall correlation $\kappa > 0$ as in Definition \ref{def:ord_opt}, and its associated Gaussian copula correlation $\rho$. It is to be kept in mind that we can take \eqref{eq:substitution} to bring the context back into controller tuning. Also, the standing assumptions can be stated in an analogous way for the distribution $\left(Z, X\right)$. In particular, the analogy to Assumption \ref{assump:sampling} is to let an i.i.d. sample of size $n$ be denoted by
\begin{multline}
\left(Z_{1}, X_{1}\right), \dots, \left(Z_{n}, X_{n}\right) \\
\underset{\mathrm{st}}{=} \left(\overline{J}\left(\theta_{1}\right), J\left(\psi_{1}, \theta_{1}\right)\right), \dots, \left(\overline{J}\left(\theta_{n}\right), J\left(\psi_{n}, \theta_{n}\right)\right).
\label{eq:bivariate_sample}
\end{multline}

First, we consider the following point estimators for $\alpha$ and $\rho$.

\begin{definition}[Point estimator for $\alpha$]
From the sample \eqref{eq:bivariate_sample}, a point estimate of $\alpha = \operatorname{Pr}\left(X \leq x^{*}\right)$ for performance threshold $x^{*}$ is
\begin{equation}
\widehat{\alpha}_{n} := \dfrac{1}{n}\sum_{i = 1}^{n}\mathbb{I}_{\left\{X_{i} \leq x^{*}\right\}}.
\label{eq:point_estimate_alpha}
\end{equation}
\end{definition}

\begin{definition}[Point estimator for $\rho$]
From the sample \eqref{eq:bivariate_sample}, a point estimate of the correlation $\rho$ for the associated Gaussian copula is
\begin{equation}
\widehat{\rho}_{n} := \sin\left(\dfrac{\pi}{2}\max\left\{0, \widehat{\kappa}_{n}\right\}\right),
\label{eq:point_estimate_rho}
\end{equation}
where $\widehat{\kappa}_{n}$ is the sample Kendall correlation
\begin{equation}
\widehat{\kappa}_{n} := \dfrac{1}{n\left(n-1\right)}\sum_{i=1}^{n}\sum_{j=1}^{n}\operatorname{sign}\left(\left(X_{i}-X_{j}\right)\left(Z_{i}-Z_{j}\right)\right).
\label{eq:point_estimate_kappa}
\end{equation}
\end{definition}

We also provide monotonicity properties of the OO success probability in $\alpha$ and $\rho$, analogous to Theorem \ref{thm:OO_properties}\ref{thm:OO_monotonicity_m}, \ref{thm:OO_monotonicity_n}.

\begin{lemma}[Monotonicity in $\alpha$]
Given the pair $\left(\bar{n}, \bar{m}\right) \in \mathbb{N} \times \left[1, n\right]$, then for all $\alpha \leq \alpha'$ such that $\alpha \in \left(0, 1\right]$ and $\alpha' \in \left(0, 1\right]$, we have for the OO success probability \eqref{eq:ord-opt-defn} that
\begin{equation}
p_{\mathrm{success}}\left(\bar{n}, \bar{m}, \alpha\right) \leq p_{\mathrm{success}}\left(\bar{n}, \bar{m}, \alpha'\right).
\end{equation}
\label{lem:monotonicity_alpha}
\end{lemma}
\begin{pf}
By De Morgan's laws (i.e. complement of the union is the intersection of the complements), put the definition of $p_{\mathrm{success}}\left(n, m, \alpha\right)$ from \eqref{eq:ord-opt-defn} in terms of
\begin{equation}
p_{\mathrm{success}}\left(n, m, \alpha\right) = 1 - \operatorname{Pr}\left(\bigcap_{i = 1}^{m}\left\{X_{\left\langle i\right\rangle} > x_{\alpha}^{*}\right\}\right).
\end{equation}
Then apply the properties that $x_{\alpha}^{*}$ is non-decreasing in $\alpha$ and $\operatorname{Pr}\left(\bigcap_{i = 1}^{m}\left\{X_{\left\langle i\right\rangle} > x_{\alpha}^{*}\right\}\right)$ is non-increasing in $x_{\alpha}^{*}$. 
\qed
\end{pf}

\begin{lemma}[Monotonicity in $\rho$]
Given the triplet $\left(\bar{n}, \bar{m}, \bar{\alpha}\right) \in \mathbb{N} \times \left[1, n\right] \times \left(0, 1\right]$, then for all $\rho \leq \rho'$ such that $\rho \in \left(0, 1\right]$ and $\rho' \in \left(0, 1\right]$, we have for the Gaussian copula OO success probability
\begin{equation}
p_{\mathrm{success}}^{\mathcal{N}}\left(\bar{n}, \bar{m}, \bar{\alpha}, \rho\right) \leq p_{\mathrm{success}}^{\mathcal{N}}\left(\bar{n}, \bar{m}, \bar{\alpha}, \rho'\right).
\end{equation}
\label{lem:monotonicity_rho}
\end{lemma}
\begin{pf}
Provided in Appendix \ref{sec:monotonicity_rho}.
\qed
\end{pf}

Confidence bounds for $\alpha$ and $\rho$ can be obtained from the following concentration inequalities.

\begin{lemma}[Concentration inequalities for $\alpha$]
For $a > 0$, we have
\begin{gather}
\operatorname{Pr}\left(\widehat{\alpha}_{n}-\alpha<-a\right) \leq \exp\left(-2na^{2}\right) \label{eq:conc_ineq_alpha_lower} \\
\operatorname{Pr}\left(\widehat{\alpha}_{n}-\alpha>a\right) \leq \exp\left(-2na^{2}\right). \label{eq:conc_ineq_alpha_upper}
\end{gather}
\label{lem:concentration_alpha}
\end{lemma}
\begin{pf}
Recognising that $n\widehat{\alpha}_{n}$ is a sum of independent Bernoulli random variables (each bounded between $0$ and $1$) with mean $\alpha$, use Hoeffding's inequality \cite[Theorem 1.1]{Dubhashi2009} to obtain
\begin{align}
\operatorname{Pr}\left(\widehat{\alpha}_{n}-\alpha>a\right) &= \operatorname{Pr}\left(n\widehat{\alpha}_{n}-n\alpha>na\right) \\
&\leq \exp\left(-2na^{2}\right),
\end{align}
and analogously for the lower tail bound. 
\qed
\end{pf}

\begin{lemma}[Concentration inequalities for $\rho$]
Under Assumption \ref{assump:copula} with the substitution \eqref{eq:substitution}, for $r > 0$, we have
\begin{gather}
\operatorname{Pr}\left(\widehat{\rho}_{n}-\rho<-r\right) \leq \exp\left(-\left\lfloor \dfrac{n}{2}\right\rfloor \dfrac{2r^{2}}{\pi^{2}}\right) \label{eq:conc_ineq_rho_lower} \\
\operatorname{Pr}\left(\widehat{\rho}_{n}-\rho>r\right) \leq \exp\left(-\left\lfloor \dfrac{n}{2}\right\rfloor \dfrac{2r^{2}}{\pi^{2}}\right). \label{eq:conc_ineq_rho_upper}
\end{gather}
\label{lem:concentration_rho}
\end{lemma}
\begin{pf}
Provided in Appendix \ref{sec:concentration_rho}.
\qed
\end{pf}
\begin{remark}
A two-tailed bound similar to Lemma \ref{lem:concentration_rho} with a slightly different exponent can be found in \cite[Theorem 4.2]{Liu2012a}. Applying the fact that $n/4 \leq \left\lfloor n/2\right\rfloor $ for all $n > 1$, one can eliminate the floor operator in \eqref{eq:conc_ineq_rho_lower}, \eqref{eq:conc_ineq_rho_upper} and recover the same exponent as found in \cite{Liu2012a}.
\end{remark}
From the upper tailed concentration inequalities for $\alpha$ and $\rho$, we may then derive lower confidence bounds. To derive a lower confidence bound for $\alpha$ with confidence level at least $1 - \beta_{1}$, equate $\exp\left(-2na^{2}\right) = \beta_{1}$ and rearrange in the upper-tailed bound \eqref{eq:conc_ineq_alpha_upper} to obtain
\begin{equation}
\operatorname{Pr}\left(\widehat{\alpha}_{n}-\alpha> \sqrt{\dfrac{\log\left(1/\beta_{1}\right)}{2n}}\right) \leq \beta_{1}.
\end{equation}
Let
\begin{equation}
b_{1} := \sqrt{\dfrac{\log\left(1/\beta_{1}\right)}{2n}},
\end{equation}
so that
\begin{equation}
\operatorname{Pr}\left(\alpha> \widehat{\alpha}_{n} - b_{1}\right) \geq 1 - \beta_{1}.
\label{eq:high_confidence_alpha}
\end{equation}
Thus the lower confidence bound for $\alpha$ with confidence at least $1 - \beta_{1}$  is obtained as
\begin{equation}
\underline{\widehat{\alpha}}_{n} := \widehat{\alpha}_{n} - b_{1}. \label{eq:lcb_alpha}
\end{equation}
To derive a lower confidence bound for $\rho$ with confidence level at least $1 - \beta_{2}$, equate $\exp\left(-\left\lfloor n/2\right\rfloor 2r^{2}/\pi^{2}\right) = \beta_{2}$ and rearrange in the upper-tailed bound \eqref{eq:conc_ineq_rho_upper} to obtain
\begin{equation}
\operatorname{Pr}\left(\widehat{\rho}_{n}-\rho> \pi\sqrt{\log\left(\dfrac{1}{\beta_{2}}\right)\cdot\dfrac{1}{2\left\lfloor \frac{n}{2} \right\rfloor}}\right) \leq \beta_{2}.
\end{equation}
Let
\begin{equation}
b_{2} := \pi\sqrt{\log\left(\dfrac{1}{\beta_{2}}\right)\cdot\frac{1}{2\left\lfloor \frac{n}{2} \right\rfloor}},
\end{equation}
so that
\begin{equation}
\operatorname{Pr}\left(\rho > \widehat{\rho}_{n} - b_{2}\right) \geq 1 - \beta_{2}.
\label{eq:high_confidence_rho}
\end{equation}
Thus the lower confidence bound for $\rho$ with confidence at least $1 - \beta_{2}$ is obtained as
\begin{equation}
\underline{\widehat{\rho}}_{n} := \widehat{\rho}_{n} - b_{2}. \label{eq:lcb_rho}
\end{equation}

We may also bound the difference in the success probability from that of its associated Gaussian copula.

\begin{lemma}[Difference in OO success probability]
Consider the OO success probability \eqref{eq:ord-opt-defn} from Definition \ref{def:ord_opt}, and let $\rho$ be the correlation of the associated Gaussian copula. If Assumption \ref{assump:conditional_KS1} holds under the substitution \eqref{eq:substitution}, then for all $n \in \mathbb{N}$ and $\alpha \in \left(0, 1\right]$ we have
\begin{equation}
p_{\mathrm{success}}^{\mathcal{N}}\left(n, 1, \alpha, \rho\right)  - p_{\mathrm{success}}\left(n, 1, \alpha\right) \leq \nu.
\label{eq:diff_success_prob}
\end{equation}
\label{lem:diff_success_prob}
\end{lemma}
\begin{pf}
Let $\left(\widetilde{Z}, \widetilde{X}\right)$ denote the copula of $\left(Z, X\right)$ and let $\left(\widetilde{Z}, \widetilde{X}'\right)$ denote the associated Gaussian copula, where the marginal $\widetilde{Z}$ can be shared since it is a uniform $\left(0, 1\right)$ random variable. Using the fact that the first order statistic of an i.i.d. uniform $\left(0, 1\right)$ sample is $\operatorname{Beta}\left(1, n\right)$ distributed \cite[\S 1.1]{Arnold2008}, and recognising that the OO success probability only depends on the underlying copula of the distribution, we have in the case $m = 1$ that
\begin{equation}
p_{\mathrm{success}}\left(n, 1, \alpha\right) = \int_{0}^{1}\operatorname{Pr}\left(\widetilde{X}\leq\alpha\middle|\widetilde{Z}=z\right)f_{U_{1:n}}\left(z\right)dz, \label{eq:psuccess_expression}
\end{equation}
where $f_{U_{1:n}}\left(\cdot\right)$ is the density of the $\operatorname{Beta}\left(1, n\right)$ distribution. Likewise
\begin{equation}
p_{\mathrm{success}}^{\mathcal{N}}\left(n, 1, \alpha, \rho\right) = \int_{0}^{1}\operatorname{Pr}\left(\widetilde{X}'\leq\alpha\middle|\widetilde{Z}=z\right)f_{U_{1:n}}\left(z\right)dz. \label{eq:psuccess_expression_gaussian}
\end{equation}
The difference between these is
\begin{align}
& p_{\mathrm{success}}^{\mathcal{N}}\left(n,1,\alpha,\rho\right)-p_{\mathrm{success}}\left(n,1,\alpha\right) \\
&= \begin{multlined}[t] \int_{0}^{1}\left(\operatorname{Pr}\left(\widetilde{X}'\leq\alpha\middle|\widetilde{Z}=z\right) \right. \\
\left. -\operatorname{Pr}\left(\widetilde{X}\leq\alpha\middle|\widetilde{Z}=z\right)\right)f_{U_{1:n}}\left(z\right)dz
\end{multlined} \label{eq:diff_success_prob1} \\
&\leq \begin{multlined}[t] \int_{0}^{1}\sup_{\alpha\in\left(0,1\right]}\left\{ \operatorname{Pr}\left(\widetilde{X}'\leq\alpha\middle|\widetilde{Z}=z\right) \right. \\
\left. -\operatorname{Pr}\left(\widetilde{X}\leq\alpha\middle|\widetilde{Z}=z\right)\right\} f_{U_{1:n}}\left(z\right)dz
\end{multlined} \label{eq:diff_success_prob2} \\
&\leq \nu\int_{0}^{1}f_{U_{1:n}}\left(z\right)dz \\
&= \nu,
\end{align}
where the second inequality is from \eqref{eq:conditional_KS1} in Assumption \ref{assump:conditional_KS1}.
\qed
\end{pf}

Using Lemma \ref{lem:diff_success_prob}, the aforementioned properties on $\alpha$ and $\rho$, as well as the lower bound for $p_{\mathrm{success}}$ in \eqref{eq:success-prob-copula-lb-optimised}, we are ready to establish a lower confidence bound on the OO success probability.

\begin{thm}[Lower confidence bound for $p_{\mathrm{success}}$]
Consider the OO success probability \eqref{eq:ord-opt-defn} from Definition \ref{def:ord_opt}. If Assumption \ref{assump:conditional_KS1} holds under the substitution \eqref{eq:substitution}, then from the sample \eqref{eq:bivariate_sample}, with confidence at least $1 - \beta_{1} - \beta_{2}$, we have
\begin{equation}
\widehat{p}_{\mathrm{success}}^{\mathcal{N}}\left(n, 1, \underline{\widehat{\alpha}}_{n}, \underline{\widehat{\rho}}_{n}\right)  - \nu \leq p_{\mathrm{success}}\left(n, m, \alpha\right).
\end{equation}
\label{thm:lcb_psuccess}
\end{thm}
\begin{pf}
As $\widehat{p}_{\mathrm{success}}^{\mathcal{N}}$ from \eqref{eq:success-prob-copula-lb-optimised} is a lower bound, then
\begin{equation}
\widehat{p}_{\mathrm{success}}^{\mathcal{N}}\left(n, 1, \underline{\widehat{\alpha}}_{n}, \underline{\widehat{\rho}}_{n}\right) \leq p_{\mathrm{success}}^{\mathcal{N}}\left(n, 1, \underline{\widehat{\alpha}}_{n}, \underline{\widehat{\rho}}_{n}\right).
\label{eq:lcb_psuccess_lb}
\end{equation}
Applying Lemma \ref{lem:diff_success_prob} (which requires Assumption \ref{assump:conditional_KS1}), this implies
\begin{equation}
\widehat{p}_{\mathrm{success}}\left(n, 1, \underline{\widehat{\alpha}}_{n}, \underline{\widehat{\rho}}_{n}\right) - \nu \leq p_{\mathrm{success}}'\left(n, 1, \underline{\widehat{\alpha}}_{n}, \underline{\widehat{\rho}}_{n}\right).
\label{eq:lcb_psucess_implication}
\end{equation}
Note that the property of monotonicity in $m$ from Theorem \ref{thm:OO_properties}\ref{thm:OO_monotonicity_m} also applies to any copula, because from \eqref{eq:ord-opt-defn}, we see that
\begin{equation}
\operatorname{Pr}\left(\bigcup_{i = 1}^{m}\left\{X_{\left\langle i\right\rangle} \leq x_{\alpha}^{*}\right\}\right) \leq \operatorname{Pr}\left(\bigcup_{i = 1}^{m + 1}\left\{X_{\left\langle i\right\rangle} \leq x_{\alpha}^{*}\right\}\right).
\end{equation}
Hence we have
\begin{equation}
p_{\mathrm{success}}\left(n, 1, \alpha\right) \leq p_{\mathrm{success}}\left(n, m, \alpha\right).
\label{eq:lcb_psuccess_monotonicity_m}
\end{equation}
Therefore
\begin{multline}
\operatorname{Pr}\left(\widehat{p}_{\mathrm{success}}^{\mathcal{N}}\left(n,1,\underline{\widehat{\alpha}}_{n},\underline{\widehat{\rho}}_{n}\right)-\nu\leq p_{\mathrm{success}}\left(n,m,\alpha\right)\right) \color{black} \\
\begin{aligned}[t] 
 & \geq\operatorname{Pr}\left(\widehat{p}_{\mathrm{success}}^{\mathcal{N}}\left(n,1,\underline{\widehat{\alpha}}_{n},\underline{\widehat{\rho}}_{n}\right)-\nu\leq p_{\mathrm{success}}\left(n,1,\alpha\right)\right)\\
 & \geq \operatorname{Pr}\left(\widehat{p}_{\mathrm{success}}^{\mathcal{N}}\left(n,1,\underline{\widehat{\alpha}}_{n},\underline{\widehat{\rho}}_{n}\right)\leq p_{\mathrm{success}}^{\mathcal{N}}\left(n,1,\alpha,\rho\right)\right) \\
 & \geq\operatorname{Pr}\left(p_{\mathrm{success}}^{\mathcal{N}}\left(n,1,\underline{\widehat{\alpha}}_{n},\underline{\widehat{\rho}}_{n}\right)\leq p_{\mathrm{success}}^{\mathcal{N}}\left(n,1,\alpha,\rho\right)\right)\\
 & \geq\operatorname{Pr}\left(\underline{\widehat{\alpha}}_{n}\leq\alpha, \underline{\widehat{\rho}}_{n}\leq\rho\right)\\
 & =1-\operatorname{Pr}\left(\underline{\widehat{\alpha}}_{n}>\alpha\text{ or }\underline{\widehat{\rho}}_{n}>\rho\right)\\
 & \geq1-\operatorname{Pr}\left(\underline{\widehat{\alpha}}_{n}>\alpha\right)-\operatorname{Pr}\left(\underline{\widehat{\rho}}_{n}>\rho\right)\\
 & \geq1-\beta_{1}-\beta_{2},
\end{aligned}
\label{eq:lcb_psuccess}
\end{multline}
where the first inequality is from applying \eqref{eq:lcb_psuccess_monotonicity_m}, the second inequality is due to the implication \eqref{eq:lcb_psucess_implication}, the third inequality is from \eqref{eq:lcb_psuccess_lb}, the fourth inequality is by applying the monotonicity properties from Lemmas \ref{lem:monotonicity_alpha} and \ref{lem:monotonicity_rho}, the fifth inequality is by the union bound (Boole's inequality), and the last inequality is from the lower confidence bound properties \eqref{eq:high_confidence_alpha}, \eqref{eq:lcb_alpha}, \eqref{eq:high_confidence_rho}, \eqref{eq:lcb_rho}. 
\qed
\end{pf}

\begin{remark}
A $1 - \beta_{1} - \beta_{2}$ lower confidence bound for the OO success probability under the associated Gaussian copula is $\widehat{p}_{\mathrm{success}}^{\mathcal{N}}\left(n,1,\underline{\widehat{\alpha}}_{n},\underline{\widehat{\rho}}_{n}\right)$, i.e.
\begin{multline}
\operatorname{Pr}\left(\widehat{p}_{\mathrm{success}}^{\mathcal{N}}\left(n,1,\underline{\widehat{\alpha}}_{n},\underline{\widehat{\rho}}_{n}\right)\leq p_{\mathrm{success}}^{\mathcal{N}}\left(n,1,\alpha,\rho\right)\right) \\
 \geq 1-\beta_{1}-\beta_{2}.
\end{multline}
\label{rem:lcb_psuccess}
\end{remark}

\section{Sequential Learning Algorithm}
\label{sec:algorithm}

In view of Remark \ref{rem:lcb_psuccess}, we present Algorithm \ref{alg:stopping}, which sequentially draws samples from $\left(Z, X\right)$ and stops after a random $\tau$ samples until an associated Gaussian copula OO success probability of at least $1 - \delta$ is reached, with confidence of at least $1 - \beta_{1} - \beta_{2}$. Note that this algorithm works irrespective of the value of $\nu$ in Assumption \ref{assump:conditional_KS1}, because the algorithm considers only the associated Gaussian copula. 
\begin{algorithm}[H]
\caption{Sequential learning for Gaussian copula OO success probability}
\begin{algorithmic}[1]
\Require $\delta \in \left(0, 1\right]$, $\beta_{1} \in \left(0, 1\right]$, $\beta_{2} \in \left(0, 1\right]$, performance threshold $x^{*}$, initial sample \eqref{eq:bivariate_sample} of size $n$
\State $n \gets n + 1$
\State Independently sample $\left(Z, X\right)$ and add to existing samples
\State Compute $\widehat{\alpha}_{n}$, $\widehat{\rho}_{n}$ via \eqref{eq:point_estimate_alpha}, \eqref{eq:point_estimate_rho}, \eqref{eq:point_estimate_kappa}
\State Compute $\underline{\widehat{\alpha}}_{n}$, $\underline{\widehat{\rho}}_{n}$ via \eqref{eq:lcb_alpha}, \eqref{eq:lcb_rho} using $\beta_{1}$, $\beta_{2}$ respectively
\State $p \gets \widehat{p}_{\mathrm{success}}^{\mathcal{N}}\left(n, 1, \underline{\widehat{\alpha}}_{n}, \underline{\widehat{\rho}}_{n}\right)$
\State If $p \geq 1 - \delta$, continue, otherwise go to step 1
\State $\tau \gets n$
\end{algorithmic}
\label{alg:stopping}
\end{algorithm}

Qualitatively, as $n$ increases, the confidence widths $b_{1}$ and $b_{2}$ decrease to zero. The lower bound from Theorem \ref{thm:OO_properties}\ref{thm:OO_success_prob_copula_lb} also stipulates that $\widehat{p}_{\mathrm{success}}^{\mathcal{N}}$ is increasing in $n$. Thus, we intuitively reason that Algorithm \ref{alg:stopping} eventually terminates with sufficiently large $n$. This intuition can be made precise with the following theorem and subsequent corollary, which uses the concentration inequalities for $\alpha$ and $\rho$ to bound the distribution of the time at which Algorithm \ref{alg:stopping} stops.

\begin{thm}[Bound on stopping time]
Fix $\delta$, $\beta_{1}$, $\beta_{2}$ in Algorithm \ref{alg:stopping}. Given some $n \in \mathbb{N}$, suppose the pair $\left(\alpha^{*}, \rho^{*}\right)$ satisfies $\widehat{p}_{\mathrm{success}}^{\mathcal{N}}\left(n, 1, \alpha^{*}, \rho^{*}\right) \geq 1 - \delta$. Also let
\begin{gather}
a := \alpha_{0} - \alpha^{*} \label{eq:stopping_time_a} \\
r := \rho_{0} - \rho^{*}, \label{eq:stopping_time_r}
\end{gather}
where $\alpha_{0}$, $\rho_{0}$ are the actual values of $\alpha$, $\rho$ respectively.  Then, for all $n$ greater than the initial sample size, we have
\begin{multline}
56 \geq 1-\exp\left(-2n\left(\alpha_{0}-\alpha^{*}-b_{1}\right)^{2}\right) \\
-\exp\left(-\left\lfloor \dfrac{n}{2}\right\rfloor \dfrac{2\left(\rho_{0}-\rho^{*}-b_{2}\right)^{2}}{\pi^{2}}\right),
\label{eq:lb_stopping}
\end{multline}
provided $\alpha_{0} -\alpha^{*}-b_{1} > 0$ and $\rho_{0}-\rho^{*}-b_{2} > 0$.
\label{thm:stopping}
\end{thm}
\begin{pf}
\begingroup
\allowdisplaybreaks
We may bound
\begin{align}
\operatorname{Pr}\left(\tau \leq n\right) &\geq \operatorname{Pr}\left(\widehat{p}_{\mathrm{success}}^{\mathcal{N}}\left(n, 1, \underline{\widehat{\alpha}}_{n}, \underline{\widehat{\rho}}_{n}\right) \geq 1-\delta\right) \\
& \geq \operatorname{Pr}\left(\underline{\widehat{\alpha}}_{n} \geq \alpha^{*}, \underline{\widehat{\rho}}_{n} \geq \rho^{*}\right) \\
&= 1-\operatorname{Pr}\left(\underline{\widehat{\alpha}}_{n}<\alpha^{*}\text{ or }\underline{\widehat{\rho}}_{n}<\rho^{*}\right) \\
&\geq \begin{multlined}[t] 1-\operatorname{Pr}\left(\underline{\widehat{\alpha}}_{n}<\alpha^{*}\right) \\
-\operatorname{Pr}\left(\underline{\widehat{\rho}}_{n}<\rho^{*}\right) \end{multlined} \\
&= \begin{multlined}[t] 1-\operatorname{Pr}\left(\underline{\widehat{\alpha}}_{n}-\alpha_{0} < -a\right) \\
-\operatorname{Pr}\left(\underline{\widehat{\rho}}_{n}-\rho_{0} <-r\right)
\end{multlined} \\
&= \begin{multlined}[t] 1 - \operatorname{Pr}\left(\widehat{\alpha}_{n}-\alpha<-a+b_{1}\right) \\
-\operatorname{Pr}\left(\widehat{\rho}_{n}-\rho<-r+b_{2}\right)
\end{multlined} \\
&\geq \begin{multlined}[t] 1 - \exp\left(-2n\left(a-b_{1}\right)^{2}\right) \\
-\exp\left(-2\left\lfloor \dfrac{n}{2}\right\rfloor \dfrac{2\left(r-b_{2}\right)^{2}}{\pi^{2}}\right),
\end{multlined}
\end{align}
where the first inequality holds because of the stopping condition, the second inequality is by definition of $\alpha^{*}$ and $\rho^{*}$ along with monotonicity properties from Lemmas \ref{lem:monotonicity_alpha} and \ref{lem:monotonicity_rho}, the third inequality is from the union bound (Boole's inequality), and the fourth equality is by application of the lower tailed concentration inequalities \eqref{eq:conc_ineq_alpha_lower}, \eqref{eq:conc_ineq_rho_lower} from Lemmas \ref{lem:concentration_alpha} and \ref{lem:concentration_rho} respectively. Substituting \eqref{eq:stopping_time_a}, \eqref{eq:stopping_time_r} completes the proof.
\endgroup
\qed
\end{pf}

\begin{corollary}[Finite stopping time]
Under Assumptions  \ref{assump:copula} and \ref{assump:feasible} with the substitution \eqref{eq:substitution}, the stopping time $\tau$ from Algorithm \ref{alg:stopping} satisfies
\begin{equation}
\operatorname{Pr}\left(\tau < \infty\right) = 1.
\end{equation}
\label{cor:finite_stopping}
\end{corollary}
\begin{pf}
Assumptions \ref{assump:copula} and \ref{assump:feasible} ensure that $\alpha_{0} > 0$ and $\rho_{0} > 0$. By Theorem \ref{thm:OO_properties}\ref{thm:OO_monotonicity_n}, \ref{thm:high_probability_success}, for any $\delta > 0$ there exists a pair $\left(\alpha^{*}, \rho^{*}\right)$ such that $\alpha_{0}-\alpha^{*}-b_{1} > 0$ and $\rho_{0}-\rho^{*}-b_{2} > 0$ for all $n$ greater than some sufficiently large number. Hence from the monotone convergence theorem \cite[Theorem 4.8]{Capinski2004}, we have
\begingroup
\allowdisplaybreaks
\begin{multline}
\operatorname{Pr}\left(\tau < \infty\right) \\
\begin{aligned}
&= \lim_{n \to \infty}\operatorname{Pr}\left(\tau<n+1\right) \\
&= \lim_{n \to \infty}\operatorname{Pr}\left(\tau\leq n\right) \\
&\geq \begin{multlined}[t] \lim_{n \to \infty}\left[ 1 - \exp\left(-2n\left(\alpha_{0}-\alpha^{*}-b_{1}\right)^{2}\right) \phantom{\dfrac{2\left(\rho_{0}\right)^{2}}{\pi^{2}}} \right. \\
\left. -\exp\left(-\left\lfloor \dfrac{n}{2}\right\rfloor \dfrac{2\left(\rho_{0}-\rho^{*}-b_{2}\right)^{2}}{\pi^{2}}\right) \right]
\end{multlined} \\
&= 1,
\end{aligned}
\end{multline}
\endgroup
where the inequality is by applying Theorem \ref{thm:stopping}. \qed
\end{pf}

\begin{remark}[Optimised bound on stopping time]
We can also numerically optimise the bound \eqref{eq:lb_stopping} with respect to $\left(\alpha^{*}, \rho^{*}\right)$. This can be useful for characterising the computational requirement (i.e. number of samples needing to be simulated) of the algorithm. Further details on optimising the bound are provided in Appendix \ref{sec:optimised_bound}.
\label{ref:optimised_bound}
\end{remark}

\subsection{Controller Tuning Algorithm}

Next, we specialise Algorithm \ref{alg:stopping} to the context of controller tuning, in order to explicitly address Problem \ref{prob:tuning}. This is presented in Algorithm \ref{alg:tuning}, which now also outputs the tuned controller $\theta_{\tau}^{*}$. 
\begin{algorithm}[H]
\caption{Probabilistically robust controller tuning}
\begin{algorithmic}[1]
\Require $\delta \in \left(0, 1\right]$, $\beta_{1} \in \left(0, 1\right]$, $\beta_{2} \in \left(0, 1\right]$, performance threshold $J^{*}$, initial sample \eqref{eq:bivariate_sample_performances} of size $n$
\State $n \gets n + 1$
\State Independently sample $\theta_{i}$ and $\psi_{i}$
\State Form $\left(Z_{i}, X_{i}\right) \gets \left(\overline{J}\left(\theta_{i}\right), J\left(\psi_{i}, \theta_{i}\right)\right)$ and add to existing samples
\State Compute $\widehat{\alpha}_{n}$ via \eqref{eq:point_estimate_alpha} with performance threshold $J^{*}$
\State Compute $\widehat{\rho}_{n}$ via \eqref{eq:point_estimate_rho}, \eqref{eq:point_estimate_kappa}
\State Compute $\underline{\widehat{\alpha}}_{n}$, $\underline{\widehat{\rho}}_{n}$ via \eqref{eq:lcb_alpha}, \eqref{eq:lcb_rho} using $\beta_{1}$, $\beta_{2}$ respectively
\State $p \gets \widehat{p}_{\mathrm{success}}^{\mathcal{N}}\left(n, 1, \underline{\widehat{\alpha}}_{n}, \underline{\widehat{\rho}}_{n}\right)$
\State If $p \geq 1 - \delta$, continue, otherwise go to step 1
\State $\tau \gets n$
\State $\theta_{\tau}^{*} \gets \argmin_{\theta_{i} \in \left\{\theta_{1}, \dots, \theta_{\tau}\right\}}\overline{J}\left(\theta_{i}\right)$
\end{algorithmic}
\label{alg:tuning}
\end{algorithm}

By chaining the confidence level with the OO success probability, we demonstrate how Algorithm \ref{alg:tuning} addresses Problem \ref{prob:tuning}, via the following theorem.

\begin{thm}
Suppose Algorithm \ref{alg:tuning} is applied to tuning controllers of a system with a performance function $J\left(\psi, \theta\right)$. Let $\theta_{\tau}^{*}$ denote the candidate solution output by the algorithm. Under Assumptions \ref{assump:copula}, \ref{assump:conditional_KS1}, \ref{assump:feasible} and \ref{assump:sampling}, then given any $\gamma \in \left(\nu, 1\right]$ and $\delta > 0$, $\beta_{1} > 0$, $\beta_{2} > 0$ with $\delta + \beta_{1} + \beta_{2} = \gamma - \nu$, then
\begin{equation}
\operatorname{Pr}_{\psi^{\sharp}, \theta_{\tau}^{*}}\left(J\left(\psi^{\sharp},\theta_{\tau}^{*}\right)\leq J^{*}\right) \geq 1 - \gamma.
\end{equation}
\label{thm:tuning}
\end{thm}
\begin{pf}
Let $\widetilde{n}_{\delta}\left(\alpha, \rho\right)$ denote the smallest integer $n$  such that $p_{\mathrm{success}}^{\mathcal{N}}\left(n, 1, \alpha, \rho\right) \geq 1 - \delta$. Combining this with Lemma \ref{lem:diff_success_prob} (requiring Assumption \ref{assump:conditional_KS1}), we have
\begin{equation}
\operatorname{Pr}_{\psi^{\sharp}, \theta_{\tau}^{*}}\left(J\left(\psi^{\sharp}, \theta_{\tau}^{*}\right) \leq J^{*} \middle| \tau \geq \widetilde{n}_{\delta}\left(\alpha_{0}, \rho_{0}\right) \right) \geq 1 - \delta - \nu.
\end{equation}
Recognise that for any $\tau$ such that $p_{\mathrm{success}}^{\mathcal{N}}\left(\tau, 1, \alpha, \rho\right) \geq 1 - \delta$, this implies
\begin{equation}
\tau \geq \widetilde{n}_{\delta}\left(\alpha, \rho\right),
\label{eq:tuning_implication}
\end{equation}
by definition of $\widetilde{n}_{\delta}\left(\alpha, \rho\right)$ and due to monotonicity in $n$ (Theorem \ref{thm:OO_properties}\ref{thm:OO_monotonicity_n}). As noted in Corollary \ref{cor:finite_stopping} (requiring Assumptions \ref{assump:copula} and \ref{assump:feasible}), the algorithm stops at time $\tau$ with probability one such that
\begin{equation}
\widehat{p}_{\mathrm{success}}^{\mathcal{N}}\left(\tau, 1, \underline{\widehat{\alpha}}_{\tau}, \underline{\widehat{\rho}}_{\tau}\right) \geq 1 - \delta.
\label{eq:tuning_stopping_condition}
\end{equation}
Then by letting $\beta = \beta_{1} + \beta_{2}$, we have
\begin{multline}
1-\beta_{1}-\beta_{2} \\
\begin{aligned} &= 1-\beta \\
&\leq \operatorname{Pr}\left(\widehat{p}_{\mathrm{success}}^{\mathcal{N}}\left(\tau, 1, \underline{\widehat{\alpha}}_{\tau},\underline{\widehat{\rho}}_{\tau}\right) \leq p_{\mathrm{success}}^{\mathcal{N}}\left(\tau, 1, \alpha_{0},\rho_{0}\right)\right) \\
&\leq \operatorname{Pr}\left(1-\delta\leq p_{\mathrm{success}}^{\mathcal{N}}\left(\tau, 1, \alpha_{0},\rho_{0}\right)\right) \\
&\leq \operatorname{Pr}\left(\tau\geq\widetilde{n}_{\delta}\left(\alpha_{0},\rho_{0}\right)\right),
\end{aligned}
\end{multline}
where the first inequality is via Remark \ref{rem:lcb_psuccess}, the second inequality is due to the stopping condition \eqref{eq:tuning_stopping_condition}, and the third inequality is due to the implication \eqref{eq:tuning_implication}. Thus
\begin{multline}
\operatorname{Pr}\left(J\left(\psi^{\sharp},\theta_{\tau}^{*}\right)\leq J^{*}\right) \\
\begin{aligned} &\geq\operatorname{Pr}\left(J\left(\psi^{\sharp},\theta_{\tau}^{*}\right)\leq J^{*},\tau\geq\widetilde{n}_{\delta}\left(\alpha_{0},\rho_{0}\right)\right) \\
&= \begin{multlined}[t] \operatorname{Pr}\left(J\left(\psi^{\sharp},\theta_{\tau}^{*}\right)\leq J^{*}\middle|\tau\geq\widetilde{n}_{\delta}\left(\alpha_{0},\rho_{0}\right)\right) \\
\times \operatorname{Pr}\left(\tau\geq\widetilde{n}_{\delta}\left(\alpha_{0},\rho_{0}\right)\right)
\end{multlined} \\
&\geq\left(1-\delta - \nu \right)\left(1-\beta\right) \\
&> 1-\delta  - \nu-\beta\\
&= 1-\gamma,
\end{aligned}
\end{multline}
because $\delta, \beta_{1}, \beta_{2}$ were chosen such that
\begin{equation}
\gamma = \delta + \beta_{1} + \beta_{2} + \nu.
\label{eq:gamma}
\end{equation}
\qed
\end{pf}

\begin{remark}
The term $\nu$ is interpreted as an upper bound on the amount of performance degradation of the algorithm in the OO success probability (appearing in the left-hand side of \eqref{eq:lcb_psuccess}), due to an unfavourable deviation of the copula of $\left(\overline{J}\left(\theta_{i}\right), J\left(\psi^{\sharp}, \theta_{i}\right)\right)$ from its associated Gaussian copula. Although the value of $\nu$ is treated as given in Theorem \ref{thm:tuning}, this might not be explicitly known a priori in practice, and instead must be assumed. However, some evidence for the amount of performance degradation can be obtain ex-post from a collected sample. This is demonstrated later in the numerical example.
\end{remark}

\begin{remark}
If $\left(\overline{J}\left(\theta_{i}\right), J\left(\psi^{\sharp}, \theta_{i}\right)\right)$ is assumed to have a Gaussian copula, we may take $\nu = 0$ as per Remark \ref{rem:conditional_KS1}, so $\operatorname{Pr}_{\psi^{\sharp}, \theta_{\tau}^{*}}\left(J\left(\psi^{\sharp},\theta_{\tau}^{*}\right)\leq J^{*}\right)$ can be made arbitrarily close to one. This is because the lower tail boundary conditional CDF for the bivariate Gaussian copula is degenerate at zero for $\rho > 0$ \cite[\S 4.3.1]{Joe2014}, so in the expression \eqref{eq:psuccess_expression_gaussian} for the OO success probability, $\lim_{n\to\infty}p_{\mathrm{success}}^{\mathcal{N}}\left(n, 1, \alpha, \rho\right) = 1$. However, there exist families of bivariate copulas where the lower tail boundary conditional CDF is non-degenerate (e.g. the bivariate Frank copula \cite[\S 4.5.1]{Joe2014}), so in the expression \eqref{eq:psuccess_expression}, generally $\lim_{n\to\infty}p_{\mathrm{success}}\left(n, 1, \alpha\right) \neq 1$. Therefore in the case $\nu > 0$, we generally cannot make $\operatorname{Pr}_{\psi^{\sharp}, \theta_{\tau}^{*}}\left(J\left(\psi^{\sharp},\theta_{\tau}^{*}\right)\leq J^{*}\right)$ arbitrarily close to one.
\end{remark}

\begin{remark}
An explicit choice of algorithm settings may, for instance, be $\delta = \beta_{1} = \beta_{2} = \left(\gamma - \nu\right)/3$, or $\delta = \left(\gamma - \nu\right)/2$, $\beta_{1} = \beta_{2} = \left(\gamma - \nu\right)/4$. The lower bound on the stopping time in Theorem \ref{thm:stopping} will generally depend on $\delta$, $\beta_{1}$, $\beta_{2}$, so by choosing an appropriate combination of $\delta$, $\beta_{1}$, $\beta_{2}$ (a practice known as risk allocation \cite{Ono2008}), the performance of the algorithm can potentially be improved. However, doing so would not be reasonable in practice since it also requires the actual values of $\alpha$ and $\rho$ to be known.
\end{remark}

\section{Numerical Example}
\label{sec:numerical}

We demonstrate our proposed approach on a numerical example, of tuning model predictive controllers (MPC) on automotive diesel engines. Typically in production, a single controller will be tuned for a fleet of vehicles. However, the exact model representing any individual engine dynamics may differ slightly. Let the linearised nominal model for the engine be
\begin{equation}
\widetilde{x}_{k+1}=A^{\left(\varrho\right)}\widetilde{x}_{k}+B^{\left(\varrho\right)}\widetilde{u}_{k},
\label{eq:nominal_model}
\end{equation}
with $\varrho \in \left\{ 1,\dots,12\right\}$ representing each of the linearisation points, which is given by the current engine operating condition. Any particular engine has some variation in its dynamics, with the nominal matrices $A^{\left(\varrho\right)} \in \mathbb{R}^{\mathsf{n} \times \mathsf{n}}$, $B^{\left(\varrho\right)} \in \mathbb{R}^{\mathsf{n} \times \mathsf{m}}$ subjected to a disturbance
\begin{gather}
A^{\left(\varrho\right)}=\overline{A}^{\left(\varrho\right)}+\Delta_{A}^{\left(\varrho\right)} \label{eq:perturbed_A}\\
B^{\left(\varrho\right)}=\overline{B}^{\left(\varrho\right)}+\Delta_{B}^{\left(\varrho\right)} \label{eq:perturbed_B},
\end{gather}
for each $\varrho \in \left\{ 1,\dots,12\right\}$, where $\Delta_{A}^{\left(\varrho\right)} \in \mathbb{R}^{\mathsf{n} \times \mathsf{n}}$, $\Delta_{A}^{\left(\varrho\right)} \in \mathbb{R}^{\mathsf{n} \times \mathsf{m}}$ are random disturbances. These disturbances model the uncertainty of dynamics for vehicles in the fleet. We apply Algorithm \ref{alg:tuning} to tuning engine controllers so that the performance (measured in terms of squared tracking error), will be robust to these variations.

\subsection{System Description}

The air-path of an automotive diesel engine can be locally represented by a reduced order linear model with $\mathsf{n}=4$ states, and $\mathsf{m}=3$ inputs \cite{Shekhar2017}. The state vector is denoted by
\begin{equation}
x = \begin{bmatrix}
\mathsf{p}_{\mathrm{im}} & \mathsf{p}_{\mathrm{em}} & W_{\mathrm{comp}} & y_{\mathrm{EGR}}
\end{bmatrix}^{\top},
\end{equation}
where $\mathsf{p}_{\mathrm{im}}$ is the engine intake manifold pressure, $\mathsf{p}_{\mathrm{em}}$ is the exhaust manifold pressure, $W_{\mathrm{comp}}$ is the compressor mass flow rate and $y_{\mathrm{EGR}}$ is the known as the exhaust gas recirculation (EGR) rate. The inputs are composed of the actuation signals
\begin{equation}
u = \begin{bmatrix}
u_{\mathrm{thr}} & u_{\mathrm{EGR}} & u_{\mathrm{VGT}}
\end{bmatrix}^{\top},
\end{equation}
where $u_{\mathrm{thr}}$ is for the throttle valve, $u_{\mathrm{EGR}}$ is for the EGR valve and $u_{\mathrm{VGT}}$ is for the variable geometry turbine (VGT) vane. In the same vein as \cite{Maass2020, Sankar2020}, the MPC is designed in the perturbed state and inputs
\begin{gather}
\widetilde{x}_{k}=x_{k}-\overline{x}^{\left(\varrho\right)} \\
\widetilde{u}_{k}=u_{k}-\overline{u}^{\left(\varrho\right)},
\end{gather}
where $\overline{x}^{\left(\varrho\right)}$, $\overline{u}^{\left(\varrho\right)}$ are the given steady-state states and inputs respectively for the current operating condition. The outputs of interest are the intake manifold pressure and EGR rate, i.e.
\begin{equation}
y = \begin{bmatrix}
\mathsf{p}_{\mathrm{im}} & y_{\mathrm{EGR}}
\end{bmatrix}^{\top} = \begin{bmatrix}1 & 0 & 0 & 0\\
0 & 0 & 0 & 1
\end{bmatrix}x = Cx.
\end{equation}
The steady-state outputs $\overline{y}^{\left(\varrho\right)} = C\overline{x}^{\left(\varrho\right)}$ for a given operating condition also act as the reference outputs for that operating condition. Thus, a time-varying output reference trajectory (induced by an engine drive-cycle and corresponding trajectory for $\varrho$) translates to regulation problem in the perturbed state. So here, the MPC law given $\varrho_{k}$ is obtained by solving the receding horizon quadratic cost problem (with prediction horizon 10):
\begin{equation}
\begin{aligned}\min_{\widetilde{u}_{k|0},\dots,\widetilde{u}_{k|9}} &\begin{multlined}[t] \left\{ \sum_{i=0}^{9}\left(\widetilde{x}_{k|i}^{\top}Q^{\left(\varrho_{k}\right)}\widetilde{x}_{k|i}+\widetilde{u}_{k|i}^{\top}R^{\left(\varrho_{k}\right)}\widetilde{u}_{k|i}\right) \right. \\
\left. \phantom{\sum_{i=0}^{9}}+\widetilde{x}_{k|10}^{\top}P^{\left(\varrho_{k}\right)}\widetilde{x}_{k|10}\right\} \end{multlined} \\
\text{subject to } & \widetilde{x}_{k|i+1}=\overline{A}^{\left(\varrho_{k}\right)}\widetilde{x}_{k|i}+\overline{B}^{\left(\varrho_{k}\right)}\widetilde{u}_{k|i},\;i=0,\dots,9\\
 & \mathsf{M}^{\left(\varrho_{k}\right)}\widetilde{x}_{k|i}\leq\mathsf{f}^{\left(\varrho_{k}\right)},\;i=1,\dots,10\\
 & \mathsf{E}^{\left(\varrho_{k}\right)}\widetilde{u}_{k|i}\leq\mathsf{h}^{\left(\varrho_{k}\right)},\;i=0,\dots,9\\
 & \left|\widetilde{u}_{k|i}-\widetilde{u}_{k|i-1}\right|\leq u_{\mathrm{slew}},\;i=0,\dots,9,
\end{aligned}
\label{eq:problem_MPC}
\end{equation}
where $Q^{\left(\varrho\right)} \succ 0$, $P^{\left(\varrho\right)} \succ 0$, $R^{\left(\varrho\right)} \succ 0$ for each $\varrho \in \left\{ 1,\dots,12\right\}$, and with constraint matrices $\mathsf{M}^{\left(\varrho\right)}$, $\mathsf{f}^{\left(\varrho\right)}$, $\mathsf{E}^{\left(\varrho\right)}$, $\mathsf{h}^{\left(\varrho\right)}$ (representing physical constraints on the signals). The variable $u_{\mathrm{slew}}$ is the input slew rate. At state $\widetilde{x}_{k}$, the optimal solution to \eqref{eq:problem_MPC} with $\widetilde{x}_{k|0} = \widetilde{x}_{k}$ is obtained as $\left(\widetilde{u}_{k|0}^{*}, \dots, \widetilde{u}_{k|9}^{*}\right)$, and the input commanded at time $k$ is $\widetilde{u}_{k} = \widetilde{u}_{k|0}^{*}$.

In this particular example, we consider the task of tracking a time-varying output reference trajectory $\mathbf{y}_{\mathrm{ref}}$, that is induced by the third section of the Urban Drive-Cycle (UDC). The uncertain plant parameter $\psi$ is given by the tuple
\begin{equation}
\psi=\left(A^{\left(1\right)},B^{\left(1\right)},\dots,A^{\left(12\right)},B^{\left(12\right)}\right).
\end{equation}
To obtain its distribution $\mathcal{P}_{\psi}$, uncertainty quantification for the diesel engine air-path has been performed using a methodology based on Gaussian processes detailed in \cite{Chin2020a}, which we assume for the purpose of this example represents the uncertainty over a fleet of vehicles. As a result, the plant uncertainty is expressed as the nominal plant
\begin{equation}
\overline{\psi}=\left(\overline{A}^{\left(1\right)},\overline{B}^{\left(1\right)},\dots,\overline{A}^{\left(12\right)},\overline{B}^{\left(12\right)}\right)
\end{equation}
plus independent zero-mean Gaussian perturbations to each of the elements of the matrices, in the same way as of \eqref{eq:perturbed_A}, \eqref{eq:perturbed_B}.

The tunable controller parameters $\theta$ are the tuple of positive definite cost matrices
\begin{equation}
\theta=\left(Q^{\left(1\right)},P^{\left(1\right)},R^{\left(1\right)},\dots,Q^{\left(12\right)},P^{\left(12\right)},R^{\left(12\right)}\right).
\label{eq:controller_tuple}
\end{equation}
Let the system performance function for our tracking problem be defined as the tracking error:
\begin{equation}
J\left(\psi,\theta\right)=\left\Vert \mathcal{Y}_{\psi,\theta}-\mathbf{y}_{\mathrm{ref}}\right\Vert _{\mathsf{F}}^{2},
\end{equation}
where $\left\Vert \cdot\right\Vert_{\mathsf{F}}^{2}$ denotes the Frobenius norm, $\mathbf{y}_{\mathrm{ref}}\in\mathbb{R}^{2\times T}$ is the reference trajectory, and $\mathcal{Y}_{\psi, \theta}\in\mathbb{R}^{2\times T}$ is the actual discrete-time closed-loop trajectory of $y_{k}$ under controller $\theta$ on plant $\psi$. The outputs have been normalised to be within the same order of magnitude. Let the performance comparison function be the tracking from using the nominal model $\overline{\psi}$ as the plant dynamics in closed-loop under controller $\theta$, i.e.
\begin{equation}
\overline{J}\left(\theta\right)=\left\Vert \mathcal{Y}_{\overline{\psi},\theta}-\mathbf{y}_{\mathrm{ref}}\right\Vert _{\mathsf{F}}^{2}.
\end{equation}
The mechanism $\mathcal{P}_{\theta}$ we use for randomly generating the $Q^{\left(\varrho\right)}$, $R^{\left(\varrho\right)}$ matrices for each $\varrho\in\left\{ 1,\dots,12\right\}$, using a comparable approach to \cite{Ira2020, Maass2020}, is given by
\begin{gather}
Q^{\left(\varrho\right) } = \mathsf{W}_{Q}^{\left(\varrho\right)}\mathsf{D}_{Q}^{\left(\varrho\right)}\left(\mathsf{W}_{Q}^{\left(\varrho\right)}\right)^{\top} \\
R^{\left(\varrho\right)} = \mathsf{W}_{R}^{\left(\varrho\right)}\mathsf{D}_{R}^{\left(\varrho\right)}\left(\mathsf{W}_{R}^{\left(\varrho\right)}\right)^{\top},
\end{gather}
where
\begin{itemize}
\item $\mathsf{W}_{Q}^{\left(\varrho\right)}$ and $\mathsf{W}_{R}^{\left(\varrho\right)}$ are uniformly random orthogonal matrices of dimensions $4\times 4$ and $3\times 3$ respectively,
\item $\mathsf{D}_{Q}^{\left(\varrho\right)}$ and $\mathsf{D}_{R}^{\left(\varrho\right)}$ are a diagonal matrices whose diagonal elements are independently $\operatorname{Exp}\left(1\right)$ and $\operatorname{Exp}\left(100\right)$ distributed respectively.
\end{itemize}
Then, $P^{\left(\varrho\right)}$ is fixed with respect to $\overline{A}^{\left(\varrho\right)}, \overline{B}^{\left(\varrho\right)}, Q^{\left(\varrho\right)}, R^{\left(\varrho\right)}$ by solving the discrete-time algebraic Riccati equation.

\subsection{Single Tuned Controller for a Fleet}
\label{sec:single}

\begin{figure}[!htb]
\includegraphics[width=0.45\textwidth]{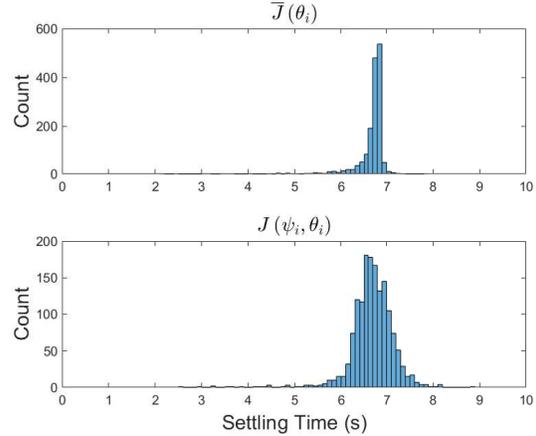}\centering
\caption{Histograms of $\overline{J}\left(\theta_{i}\right)$ and $J\left(\psi_{i}, \theta_{i}\right)$ for the 2658 samples in a single run of Algorithm \ref{alg:tuning}.}
\label{fig:histogram}
\end{figure}

We set a nominal performance threshold of $J^{*} = 1000$ to represent `good' tracking for the purposes of this example. This value was chosen before running the algorithm. After running Algorithm \ref{alg:tuning} with settings $\delta = 0.025, \beta_{1} = 0.0125, \beta_{2} = 0.0125$, then from Theorem \ref{thm:tuning} the prescribed probability of the nominal performance threshold being met in a single test is at least $0.95 - \nu$, where $\nu$ can be optimistically assumed to be small (this is validated further on). This run of the algorithm stopped after $\tau = 2658$ samples. A histogram for the performances $\overline{J}\left(\theta_{i}\right)$ and $J\left(\psi_{i}, \theta_{i}\right)$ are plotted in Figure \ref{fig:histogram}.

The best controller $\theta_{\tau}^{*}$ when evaluated on the performance comparison function $\overline{J}\left(\theta\right)$ was found to be $\overline{J}\left(\theta_{\tau}^{*}\right) =  626.06$. Upon simulating a test of this tuned controller using another randomly generated plant $\psi^{\sharp}$, we obtained a performance of $645.55$, which far outperforms the nominal threshold of $1000$.

To investigate the tuned controller performance on a fleet of vehicles, we simulated 10000 tests on another set of independently generated plants, with the same tuned controller. By the linearity of expectation, Theorem \ref{thm:tuning} prescribes that the expected proportion of tests which outperform $J^{*} = 1000$ to be at least $0.95 - \nu$. We found that all 10000 of the tests outperformed the nominal threshold, so the proportion far exceeds $0.95$. Moreover, the minimum, mean and maximum performances were $400.12$, $507.82$ and $719.01$ respectively.

\subsection{Multiple Algorithm Runs}
\label{sec:multiple_runs}

Aggregate results were also obtained for 1250 independent runs of Algorithm \ref{alg:tuning} with identical tuning procedure and settings as described above, producing 1250 tuned controllers. Each controller was then tested on its own randomly generated plant. It was found that all 1250 tests succeeded in outperforming the nominal threshold of $J^{*} = 1000$, with minimum, mean and maximum performances of $489.36$, $645.11$ and $819.42$ respectively. Note that the distribution of these 1250 tests is different from that of the 10000 tests in previous section, as each test here consists of a different controller.

\subsection{Discussion}

\begin{figure}[!htb]
\includegraphics[width=0.45\textwidth]{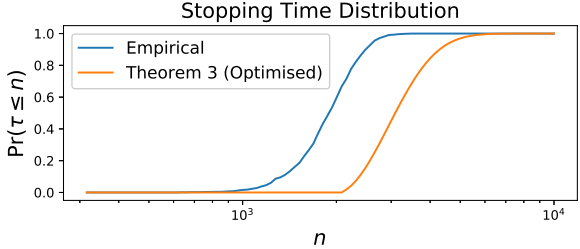}\centering
\caption{The empirical algorithm run times compared against those as suggested by Theorem \ref{thm:stopping}.}
\label{fig:stoppping}
\end{figure}

For this example, we validate Theorem \ref{thm:stopping} by plotting in Figure \ref{fig:stoppping} the numerically optimised lower bound for the CDF of the stopping time $\tau$ (using the point estimates $\widehat{\alpha}_{\tau} = 0.0700$, $\widehat{\rho}_{\tau} = 0.9792$ from the sample in Figure \ref{fig:histogram}, in place of the actual $\alpha_{0}$, $\rho_{0}$), against the empirical CDF of the stopping time for the 1250 runs. As the curves in Figure \ref{fig:stoppping} are within less than order of magnitude on the horizontal scale, this hints that Theorem \ref{thm:stopping} is not overly conservative.

To ascertain some idea of the performance degradation in $p_{\mathrm{success}}$ (i.e. left-hand side of \eqref{eq:diff_success_prob}), we may use the empirical sample from the distribution of $\left(\overline{J}\left(\theta_{i}\right), J\left(\psi_{i}, \theta_{i}\right)\right)$. A relatively large sample with size roughly $2.7 \times 10^{6}$ was compiled, taken across the multiple algorithm runs from Section \ref{sec:multiple_runs}. Then $p_{\mathrm{success}}$ was Monte-Carlo estimated according to Definition \ref{def:ord_opt} with $n = 2658$, $\alpha = 0.0700$ (i.e. using the empirical values obtained from Section \ref{sec:single}), by \textit{bootstrapping} (re-sampling) from the large sample. This was compared against $p_{\mathrm{success}}^{\mathcal{N}}$ for the associated Gaussian copula, with $\rho = 0.9792$. Both success probabilities were computed to be very close to one. Thus, there appears virtually no or very little performance degradation that has arisen from the distribution of $\left(\overline{J}\left(\theta_{i}\right), J\left(\psi^{\sharp}, \theta_{i}\right)\right)$ not having a Gaussian copula. This is consistent with having assumed a small value for $\nu$ earlier.

We also observe that the success rate of the algorithm (also close to one, as seen from Section \ref{sec:multiple_runs}) is much better than prescribed rate of up to $0.95$ (when taking $\nu = 0$). Our explanation of this phenomenon is that it can be attributed to a combination of: 1) the distribution of $\left(\overline{J}\left(\theta_{i}\right), J\left(\psi^{\sharp}, \theta_{i}\right)\right)$ is actually more favourable to the success probability than its associated Gaussian copula; and 2) there is some conservatism in the lower confidence bounds for $\alpha$, $\rho$, and in the lower bound for $p_{\mathrm{success}}^{\mathcal{N}}$.

Finally, we compare the sample complexity of our algorithm to those required by the scenario approach \cite{Calafiore2006} for solving \eqref{eq:chance_constrained_program}. Although a direct comparison is not possible (since the assumptions, use-cases, and method of operation underpinning each algorithm differ), we still illustrate that our algorithm can result in a comparable order of magnitude for the number of samples. For the present example, the tuple \eqref{eq:controller_tuple} for $\theta$ can be parametrised as a vector of dimension
\begin{equation}
d=12\left[\dfrac{\mathsf{n}\left(\mathsf{n}+1\right)}{2}+\dfrac{\mathsf{m}\left(\mathsf{m}+1\right)}{2}\right]=192,
\end{equation}
where we do not count each of the $P^{\left(\varrho\right)}$ because they are fixed with respect to the other parameters. Choosing $\varepsilon = 0.0499$ and $\nu = 0.0001$ (defined via \eqref{eq:chance_constrained_program} and \eqref{eq:feasible_set_prob} respectively) so that $1 - \varepsilon - \nu = 0.95$, we plug these values into the bound \cite[Eq. (3)]{Calafiore2006} to obtain
\begin{equation}
\left\lceil \dfrac{2}{\varepsilon}\log\dfrac{1}{\eta}+2d+\dfrac{2d}{\varepsilon}\log\dfrac{2}{\varepsilon}\right\rceil =29156
\end{equation}
for the number of plants $\psi_{i}$ to generate. Hence, $29156$ evaluations of $J\left(\cdot, \cdot\right)$ are required per verification of constraints in the resulting convex scenario program. In contrast, Figure \ref{fig:stoppping} shows that our example typically stops between $2000$ and $3000$ samples of $\left(\theta_{i}, \psi_{i}\right)$ per run of the algorithm. A smaller number of samples is beneficial, especially in our case where each evaluation of $\overline{J}\left(\cdot\right)$ and $J\left(\cdot, \cdot\right)$ is the result of a lengthy simulation.

Note that our lower bound for the stopping time \eqref{eq:lb_stopping} does not explicitly depend on the dimension $d$ of the controller parameter. However, it is more sensitive to the correlation $\rho$ which is induced by the mechanism $\mathcal{P}_{\theta}$; increasing $\rho$ naturally causes the algorithm to stop sooner. Experimenting with this lower bound, we find $29156$ to be an upper bound for the median number of samples, i.e. $\operatorname{Pr}\left(\tau\leq29156\right)\geq0.5$, when $\rho \geq 0.8011$. This suggests that our algorithm reaches a comparable number of samples to the scenario approach for the same level of probability specification, provided that there is high correlation between $\overline{J}\left(\theta_{i}\right)$ and $J\left(\psi^{\sharp}, \theta_{i}\right)$.

\section{Conclusion}
\label{sec:conclusion}

In this paper, we addressed a robust control design problem using a sequential learning algorithm, which finds (with high probability) a candidate solution that in effect satisfies the internal performance constraint of a chance-constrained program that has black-box objective and constraint functions. Our results were enabled by exploiting the statistical correlation in the sampling of $\overline{J}\left(\theta_{i}\right)$ and $J\left(\psi^{\sharp}, \theta_{i}\right)$. The algorithm was illustrated on a numerical example involving the tuning of MPC for automotive diesel engines, and showed that the tracking performance of tuned controllers was robust to plant uncertainty in both a multi-plant setting over a fleet of vehicles with a single algorithm run, and a multi-controller setting over many algorithm runs.

Our algorithm can also potentially be applied to robust stability problems. For instance, by considering linear systems for simplicity, we could set $J\left(\psi,\theta\right)\leq J^{*}$ to be the equivalent to the condition that the closed-loop system matrix is Hurwitz. Another avenue for future research is to investigate the role that the distribution $\mathcal{P}_{\theta}$ (for sampling candidate controllers) plays in the tuned controller performance. In our formulation, $\mathcal{P}_{\theta}$ is an arbitrary choice left to the practitioner. In Section \ref{sec:numerical}, $\mathcal{P}_{\theta}$ was chosen by explicitly constructing a distribution, however another option would have been to let $\theta_{i}$ be the solution output by running a randomised optimisation algorithm with objective $\overline{J}\left(\theta\right)$. Modifying $\mathcal{P}_{\theta}$ affects the distribution of $\left(\overline{J}\left(\theta_{i}\right), J\left(\psi^{\sharp}, \theta_{i}\right)\right)$, and consequently the value of $\nu$. Hence it is perhaps worthwhile to find guiding principles in designing $\mathcal{P}_{\theta}$ which will lead to more favourable probabilistic performance specifications, or alternatively, reduced computational requirements for a fixed performance specification.

\section{Acknowledgements}                               
The first author acknowledges the support of the Elizabeth \& Vernon Puzey Scholarship. Computations described in this paper were performed using the University of Birmingham's
BEAR \cite{Bear} and The University of Melbourne's Spartan \cite{Lafayette2016} high performance computing services. 

\appendix
\section{Proof of Lemma \ref{lem:monotonicity_rho}}   
\label{sec:monotonicity_rho}
\begingroup
\allowdisplaybreaks
Let $G$ be the random variable for the number of $X$-values less than or equal to the threshold $x_{\alpha}^{*}$. Conditional on $G = g$, we can write the order statistics of the $X$-values as
\begin{equation}
X_{1:n} \leq \dots \leq X_{g:n} \leq x_{\alpha}^{*} < X_{\left(g + 1\right):n} \leq \dots \leq X_{n:n}.
\end{equation}
Using the additive Gaussian noise representation of the Gaussian copula \cite[\S A.1]{Chin2021}, we can assume without loss of generality that each $X \sim \mathcal{N}\left(0, 1\right)$, and $Z$ is formed by
\begin{equation}
Z = X + Y,
\end{equation}
where $Y \sim \mathcal{N}\left(0, \xi^{2}\right)$ is independent with $X$, and $\xi^{2} = 1/\rho^{2} - 1$. Introduce the following indexing of the $Z$-values according to the ordering of the $X$-values. We denote
\begin{equation}
Z_{\left\{i\right\}} := X_{i:n} + Y_{i}
\label{eq:alternate_Z_indexing}
\end{equation}
where the $Y_{i}$ are i.i.d. $\mathcal{N}\left(0, \xi^{2}\right)$, since the $Y$-values are independent of the ordering of the $X$-values. Let $p_{\mathrm{success}|g}^{\mathcal{N}}\left(n, m, \alpha, \rho\right)$ denote the conditional success probability, given $G = g$. An equivalent characterisation of the conditional success probability can be obtained from \cite[Equation (2.19)]{Ho2007}. This way, we may write
\begin{multline}
p_{\mathrm{success}|g}^{\mathcal{N}}\left(n, m, \alpha, \rho\right) \\
= \operatorname{Pr}\left(\min\left\{ Z_{\left\{ 1\right\} },\dots,Z_{\left\{ g\right\} }\right\} \phantom{\overset{\left(m\right)}{\min}} \right. \\
\leq \left. \overset{\left(m\right)}{\min}\left\{ Z_{\left\{ g+1\right\} },\dots,Z_{\left\{ n\right\} }\right\} \middle|G=g\right),
\end{multline}
where $\overset{\left(m\right)}{\min}\left\{\cdot\right\}$ denotes the $m$\textsuperscript{th} smallest value of its arguments. Putting \eqref{eq:alternate_Z_indexing}, we have
\begingroup
\begin{multline}
p_{\mathrm{success}|g}^{\mathcal{N}}\left(n, m, \alpha, \rho\right) 
\\ \begin{aligned}[t] &= \begin{multlined}[t] \operatorname{Pr}\left(\min\left\{ Y_{1}+X_{1:n},\dots,Y_{g}+X_{g:n}\right\} \phantom{\overset{\left(m\right)}{\min}} \right. \\
\left. \leq\overset{\left(m\right)}{\min}\left\{ Y_{g+1}+X_{\left(g+1\right):n},\dots,Y_{n}+X_{n:n}\right\} \middle|G=g\right)
\end{multlined} \\
&= \begin{multlined}[t] \operatorname{Pr}\left(\min\left\{ Y_{1}+\Delta X_{1:n},\dots,Y_{g}+\Delta X_{g:n}\right\} \phantom{\overset{\left(m\right)}{\min}} \right. \\
\leq\overset{\left(m\right)}{\min}\left\{ Y_{g+1}+\Delta X_{\left(g+1\right):n}, \right. \\
\left. \phantom{\overset{\left(m\right)}{\min}} \left.\dots,Y_{n}+\Delta X_{n:n}\right\} \middle|G=g\right),
\end{multlined}
\end{aligned} \nonumber
\end{multline}
\endgroup
where $\Delta X_{i:n} := X_{i:n} - x_{\alpha}^{*}$. Now let $\widetilde{Y}_{i} \sim \mathcal{N}\left(0, 1\right)$ represent a standardised random variable, so that $Y_{i} \underset{\mathrm{st}}{=} \xi\widetilde{Y}_{i}$, and
\begingroup
\begin{align}
&p_{\mathrm{success}|g}^{\mathcal{N}}\left(n,m,\alpha,\rho\right) \\
&= \begin{multlined}[t] \operatorname{Pr}\left(\min\left\{ \xi\widetilde{Y}_{1}+\Delta X_{1:n},\dots, \xi\widetilde{Y}_{g}+\Delta X_{g:n}\right\} \phantom{\overset{\left(m\right)}{\min}} \right. \\
\leq\overset{\left(m\right)}{\min}\left\{ \xi\widetilde{Y}_{g+1}+\Delta X_{\left(g+1\right):n}, \right. \\
\left. \phantom{\overset{\left(m\right)}{\min}} \left. \dots, \xi\widetilde{Y}_{n}+\Delta X_{n:n}\right\} \middle|G=g\right)
\end{multlined} \\
&= \begin{multlined}[t] \operatorname{Pr}\left(\min\left\{ \widetilde{Y}_{1}+\dfrac{\Delta X_{1:n}}{\xi},\dots,\widetilde{Y}_{g}+\dfrac{\Delta X_{g:n}}{\xi}\right\}  \right. \\
\leq\overset{\left(m\right)}{\min}\left\{ \widetilde{Y}_{g+1}+\dfrac{\Delta X_{\left(g+1\right):n}}{\xi}, \right. \\
\left. \phantom{\overset{\left(m\right)}{\min}} \left. \dots,\widetilde{Y}_{n}+\dfrac{\Delta X_{n:n}}{\xi}\right\} \middle|G=g\right),
\end{multlined}
\end{align}
\endgroup
because $\xi > 0$. Let any fixed realisation of the random variables $\widetilde{Y}_{1}, \dots, \widetilde{Y}_{n}, \Delta X_{1:n}, \dots, \Delta X_{n:n}$ be denoted as $\widetilde{y}_{1}, \dots, \widetilde{y}_{n}, \Delta x_{1:n}, \dots, \Delta x_{n:n}$ respectively. Observe $\Delta X_{i:n} \leq 0$ for all $i \leq g$, and $\Delta X_{i:n} > 0$  for all $i > g$. So for any $\xi' < \xi$, we have
\begin{gather}
\widetilde{y}_{i}+\dfrac{\Delta x_{i:n}}{\xi'} < \widetilde{y}_{i}+\dfrac{\Delta x_{i:n}}{\xi}, \quad \forall i \leq g \\
\widetilde{y}_{i}+\dfrac{\Delta x_{i:n}}{\xi'} > \widetilde{y}_{i}+\dfrac{\Delta x_{i:n}}{\xi}, \quad \forall i > g .
\end{gather}
Therefore it follows that
\begin{multline}
\min\left\{ \widetilde{y}_{1}+\dfrac{\Delta x_{1:n}}{\xi'},\dots,\widetilde{y}_{g}+\dfrac{\Delta x_{g:n}}{\xi'}\right\} \\
 \leq\min\left\{ \widetilde{y}_{1}+\dfrac{\Delta x_{1:n}}{\xi},\dots,\widetilde{y}_{g}+\dfrac{\Delta x_{g:n}}{\xi}\right\} \label{eq:monotonicity_rho_min}
\end{multline}
\begingroup
\begin{multline}
\overset{\left(m\right)}{\min}\left\{ \widetilde{y}_{g+1}+\dfrac{\Delta x_{\left(g+1\right):n}}{\xi},\dots,\widetilde{y}_{n}+\dfrac{\Delta x_{n:n}}{\xi}\right\} \\
\leq \overset{\left(m\right)}{\min}\left\{ \widetilde{y}_{g+1}+\dfrac{\Delta x_{\left(g+1\right):n}}{\xi'},\dots,\widetilde{y}_{n}+\dfrac{\Delta x_{n:n}}{\xi'}\right\}. \label{eq:monotonicity_rho_min_m}
\end{multline}
\endgroup
Denote $\rho' = \left(\xi'^{2} + 1\right)^{-1/2}$, so that $\rho < \rho'$. Then
\begin{align}
&p_{\mathrm{success}|g}^{\mathcal{N}}\left(n,m,\alpha,\rho'\right) \nonumber \\
&= p_{\mathrm{success}|g}^{\mathcal{N}}\left(n,m,\alpha,\left(\xi'^{2}+1\right)^{-1/2}\right) \nonumber \\
&= \begin{multlined}[t] \operatorname{Pr}\left(\min\left\{ \widetilde{Y}_{1}+\dfrac{\Delta X_{1:n}}{\xi'},\dots,\widetilde{Y}_{g}+\dfrac{\Delta X_{g:n}}{\xi'}\right\} \right. \\
\leq\overset{\left(m\right)}{\min}\left\{ \widetilde{Y}_{g+1}+\dfrac{\Delta X_{\left(g+1\right):n}}{\xi'}, \right. \\
\left. \left. \dots,\widetilde{Y}_{n}+\dfrac{\Delta X_{n:n}}{\xi'}\right\} \middle|G=g\right)
\end{multlined} \nonumber \\
&\geq \begin{multlined}[t] \operatorname{Pr}\left(\min\left\{ \widetilde{Y}_{1}+\dfrac{\Delta X_{1:n}}{\xi},\dots,\widetilde{Y}_{g}+\dfrac{\Delta X_{g:n}}{\xi}\right\} \right. \\
\leq\overset{\left(m\right)}{\min}\left\{ \widetilde{Y}_{g+1}+\dfrac{\Delta X_{\left(g+1\right):n}}{\xi}, \right. \\
\left. \left. \dots,\widetilde{Y}_{n}+\dfrac{\Delta X_{n:n}}{\xi}\right\} \middle|G=g\right) 
\end{multlined} \nonumber \\
&= p_{\mathrm{success}|g}^{\mathcal{N}}\left(n,m,\alpha,\rho\right), \nonumber
\end{align}
where the inequality is from applying \eqref{eq:monotonicity_rho_min} and \eqref{eq:monotonicity_rho_min_m}. The random variable $G$ is binomial distributed with parameters $n$, $\alpha$ (i.e. not affected by the value of $\rho$), thus
\begin{multline}
p_{\mathrm{success}}^{\mathcal{N}}\left(\bar{n}, \bar{m}, \bar{\alpha}, \rho'\right) \\
\begin{aligned} &= \sum_{g=0}^{n}p_{\mathrm{success}|g}^{\mathcal{N}}\left(\bar{n}, \bar{m}, \bar{\alpha},\rho'\right)\operatorname{Pr}\left(G=g\right) \\
&\geq \sum_{g=0}^{n}p_{\mathrm{success}|g}^{\mathcal{N}}\left(\bar{n}, \bar{m}, \bar{\alpha},\rho\right)\operatorname{Pr}\left(G=g\right) \\
&= p_{\mathrm{success}}^{\mathcal{N}}\left(\bar{n}, \bar{m}, \bar{\alpha},\rho\right).
\end{aligned}
\end{multline} 

\section{Proof of Lemma \ref{lem:concentration_rho}}       
\label{sec:concentration_rho}

We prove the upper tail concentration bound; the steps for the lower tail are similar and analogous.  From Definition \ref{def:associated_gaussian_copula}, the population Kendall correlation $\kappa$ and the associated Gaussian copula correlation $\rho$ are related by $\rho = \sin\left(\pi\kappa/2\right)$. So for $r > 0$ we have
\begin{multline}
\operatorname{Pr}\left(\widehat{\rho}_{n} - \rho > r\right) \\
\begin{aligned}[t] &= \operatorname{Pr}\left(\sin\left(\dfrac{\pi}{2}\max\left\{\widehat{\kappa}_{n}, 0\right\}\right)-\sin\left(\dfrac{\pi}{2}\kappa\right)> r\right) \nonumber \\
&= \operatorname{Pr}\left(\sin\left(\dfrac{\pi}{2}\widehat{\kappa}_{n}\right)-\sin\left(\dfrac{\pi}{2}\kappa\right)> r\right).
\end{aligned}
\end{multline}
where we able to take $\max\left\{\widehat{\kappa}_{n}, 0\right\} = \widehat{\kappa}_{n}$ since $\widehat{\kappa}_{n} \geq 0$ is necessary for $\widehat{\rho}_{n} - \rho$, as $\rho > 0$ by Assumption \ref{assump:copula}. Note that the event $\frac{\pi}{2}\left(\widehat{\kappa}_{n} - \kappa\right) > r$ together with $r > 0$ implies that $\widehat{\kappa}_{n} - \kappa > 0$. Since $\sin\left(\cdot\right)$ is $1$-Lipschitz continuous, then generally
\begin{equation}
\left|\sin\left(\dfrac{\pi}{2}\widehat{\kappa}_{n}\right)-\sin\left(\dfrac{\pi}{2}\kappa\right)\right|\leq\left|\dfrac{\pi}{2}\widehat{\kappa}_{n}-\dfrac{\pi}{2}\kappa\right|.
\end{equation}
However as we have established the sign of $\widehat{\kappa}_{n} - \kappa$, then the event $\frac{\pi}{2}\left(\widehat{\kappa}_{n} - \kappa\right) > r$ together with $r > 0$ further implies that
\begin{equation}
\sin\left(\dfrac{\pi}{2}\widehat{\kappa}_{n}\right)-\sin\left(\dfrac{\pi}{2}\kappa\right) \leq \dfrac{\pi}{2}\left(\widehat{\kappa}_{n}-\kappa\right).
\end{equation}
Thus
\begin{align}
\operatorname{Pr}\left(\widehat{\kappa}_{n}-\kappa>\dfrac{2r}{\pi}\right) &= \operatorname{Pr}\left(\dfrac{\pi}{2}\left(\widehat{\kappa}_{n}-\kappa\right)>r\right) \nonumber \\
&\geq \operatorname{Pr}\left(\sin\left(\dfrac{\pi}{2}\widehat{\kappa}_{n}\right)-\sin\left(\dfrac{\pi}{2}\kappa\right)>r\right) \nonumber \\
&= \operatorname{Pr}\left(\widehat{\rho}_{n}-\rho>r\right). \label{eq:conc_ineq_rho_tau}
\end{align}
Using the fact that $\widehat{\kappa}_{n}$ is an unbiased estimator for $\kappa$, and moreover a U-statistic with a second-order kernel bounded between $-1$ and $1$, we use \cite[Equation (5.7)]{Hoeffding1963} to obtain
\begin{equation}
\operatorname{Pr}\left(\widehat{\kappa}_{n}-\kappa>\dfrac{2r}{\pi}\right)\leq\exp\left(-\left\lfloor \dfrac{n}{2}\right\rfloor \dfrac{2r^{2}}{\pi^{2}}\right).
\end{equation}
Combining with \eqref{eq:conc_ineq_rho_tau} completes our proof.
\section{Optimised Bound for Theorem \ref{thm:stopping}}
\label{sec:optimised_bound}

The lower bound \eqref{eq:lb_stopping} for the distribution of the stopping time can be optimised by
\begin{multline}
\operatorname{Pr}(\tau \leq n) \geq 1 - \min_{\left(\alpha^{*}, \rho^{*}\right) \in \mathbb{A}}\left\{\exp\left(-2n\left(\alpha_{0}-\alpha^{*}-b_{1}\right)^{2}\right) \phantom{\exp\left(-\left\lfloor \dfrac{n}{2}\right\rfloor \dfrac{2\left(\rho-\rho^{*}-b_{2}\right)^{2}}{\pi^{2}}\right)} \right. \\
\left. + \exp\left(-\left\lfloor \dfrac{n}{2}\right\rfloor \dfrac{2\left(\rho_{0} - \rho^{*}-b_{2}\right)^{2}}{\pi^{2}}\right)\right\},
\label{eq:lb_stopping_optimised0}
\end{multline}
where
\begin{multline}
\mathbb{A} := \left\{\left(\alpha^{*}, \rho^{*}\right) \in \left(0, 1\right]^{2}: \alpha_{0} - b_{1} > \alpha^{*}, \right. \\
 \rho_{0} - b_{2} > \rho^{*},  \\
\left. \phantom{\left(0, 1\right]^{2}} \widehat{p}_{\mathrm{success}}^{\mathcal{N}}\left(n, 1, \alpha^{*}, \rho^{*}\right) \geq 1- \delta \right\}
\end{multline}
is the region of $\left(0, 1\right]^{2}$ where the gaps $\left(\alpha_{0}-b_{1}\right)-\alpha^{*}$ and $\left(\rho_{0}-b_{2}\right) -\rho^{*}$ are positive. Moreover, because the bound is improving with the gaps $\left(\alpha_{0}-b_{1}\right)-\alpha^{*}$ and $\left(\rho_{0}-b_{2}\right) -\rho^{*}$, and also because of the monotonicity properties in Lemmas \ref{lem:monotonicity_alpha} and \ref{lem:monotonicity_rho}, the optimum will lie on the Pareto front $\widehat{p}_{\mathrm{success}}^{\mathcal{N}}\left(n, 1, \alpha^{*}, \rho^{*}\right) = 1- \delta$ for $\left(\alpha^{*}, \rho^{*}\right) \in \left(0, 1\right]^{2}$. For a given $\omega$, we can instead analytically determine the Pareto front along $\widehat{p}_{\mathrm{success}, \omega}^{\mathcal{N}}\left(n, 1, \alpha^{*}, \rho^{*}\right) = 1- \delta$ using \eqref{eq:success-prob-copula-lb}. Letting
\begin{equation}
\Phi\left(\dfrac{\Phi^{-1}\left(\alpha\right) - \rho\mu_{n}}{\sqrt{1 - \rho^{2} + \rho^{2}\sigma_{n}^{2}}}\right) = 1- \delta,
\end{equation}
and putting the definitions of $\mu_{n}$, $\sigma_{n}^{2}$ and rearranging, we obtain a quadratic form in $\alpha$, $\rho$ given by
\begin{equation}
d_{1}\rho^{2}+d_{2}\Phi^{-1}\left(\alpha\right)\rho+d_{3}\Phi^{-1}\left(\alpha\right)^{2}+d_{4}=0,
\end{equation}
where
\begingroup
\allowdisplaybreaks
\begin{align}
d_{1} &= \begin{multlined}[t] -\dfrac{2\log\left(nc_{1}\right)^{2}}{\log\log2}+2\log\left(nc_{1}\right) \\
-\dfrac{2c_{2}\Phi^{-1}\left(1-\delta\right)^{2}\log\left(nc_{1}\right)}{\log\log2} \\ 
+2c_{2}\Phi^{-1}\left(1-\delta\right)^{2}-\Phi^{-1}\left(1-\delta\right)^{2} \end{multlined} \\
d_{2} &= -\dfrac{4\sqrt{c_{2}}\log\left(nc_{1}\right)^{3/2}}{\log\log2}+4\sqrt{c_{2}}\log\left(nc_{1}\right)^{1/2} \\
d_{3} &= -\dfrac{2c_{2}\log\left(nc_{1}\right)}{\log\log2}+2c_{2} \\
d_{4} &= \dfrac{2c_{2}\Phi^{-1}\left(1-\delta\right)^{2}\log\left(nc_{1}\right)}{\log\log2}-2c_{2}\Phi^{-1}\left(1-\delta\right)^{2}.
\end{align}
\endgroup
Thus given $\omega$, $n$, $\delta$, we can solve for $\rho$ in terms of $\alpha$ with
\begin{multline}
\rho_{\omega}\left(\alpha\right)=\dfrac{1}{2d_{1}}\left[-\left(d_{2}\Phi^{-1}\left(\alpha\right)\right)+ \phantom{\sqrt{\left(d_{2}\Phi^{-1}\left(\alpha\right)\right)^{2}-4d_{1}\left(d_{3}\Phi^{-1}\left(\alpha\right)^{2}+d_{4}\right)}} \right. \\ 
\left. \sqrt{\left(d_{2}\Phi^{-1}\left(\alpha\right)\right)^{2}-4d_{1}\left(d_{3}\Phi^{-1}\left(\alpha\right)^{2}+d_{4}\right)}\right].
\end{multline}
Alternatively given $\omega$, $n$, $\delta$, we can solve for $\alpha$ in terms of $\rho$ with
\begin{equation}
\alpha_{\omega}\left(\rho\right)=\Phi\left(\dfrac{-\left(d_{2}\rho\right)+\sqrt{\left(d_{2}\rho\right)^{2}-4d_{3}\left(d_{1}\rho^{2}+d_{4}\right)}}{2d_{3}}\right),
\end{equation}
where we have taken the positive solutions of the quadratics since $\alpha > 0$, $\rho > 0$. We may then proceed to optimise with respect to $\alpha^{*}$ (and $\rho^{*}$ implicitly in terms of $\alpha^{*}$) with an inner minimisation for a given $\omega$, and then optimise with respect to $\omega$ in an outer minimisation. Explicitly, \eqref{eq:lb_stopping_optimised0} becomes
\begin{multline}
\operatorname{Pr}(\tau\leq n) \\
\geq1-\inf_{\omega\in\Omega_{n}}\left\{ \min_{\alpha^{*}\in\mathbb{A}'_{\omega}}\left\{ \exp\left(-2n\left(\alpha_{0}-\alpha^{*}-b_{1}\right)^{2}\right) \phantom{\exp\left(-\left\lfloor \dfrac{n}{2}\right\rfloor \dfrac{2\left(\rho_{0}-\rho^{*}-b_{2}\right)^{2}}{\pi^{2}}\right)} \right.\right. \\
\left.\left. +\exp\left(-\left\lfloor \dfrac{n}{2}\right\rfloor \dfrac{2\left(\rho_{0}-\rho_{\omega}\left(\alpha^{*}\right)-b_{2}\right)^{2}}{\pi^{2}}\right)\right\} \right\},
\label{eq:lb_stopping_optimised}
\end{multline}
where
\begin{equation}
\mathbb{A}'_{\omega} := \left\{\alpha \in \left(0, 1\right]: \alpha_{\omega}\left(\rho_{0} - b_{2}\right)  \leq \alpha \leq \alpha_{0} - b_{1} \right\},
\end{equation}
and $\omega\in\Omega_{n} \subset \left(0, \pi/2\right)$ is defined the same as in \eqref{eq:success-prob-copula-lb-optimised}. The inner minimisation in \eqref{eq:lb_stopping_optimised} is quasiconvex, thus the optimised bound is not too difficult to numerically implement. This is further illustrated in Figure \ref{fig:optimised}.
\begin{figure}[!htb]
\includegraphics[width=0.45\textwidth]{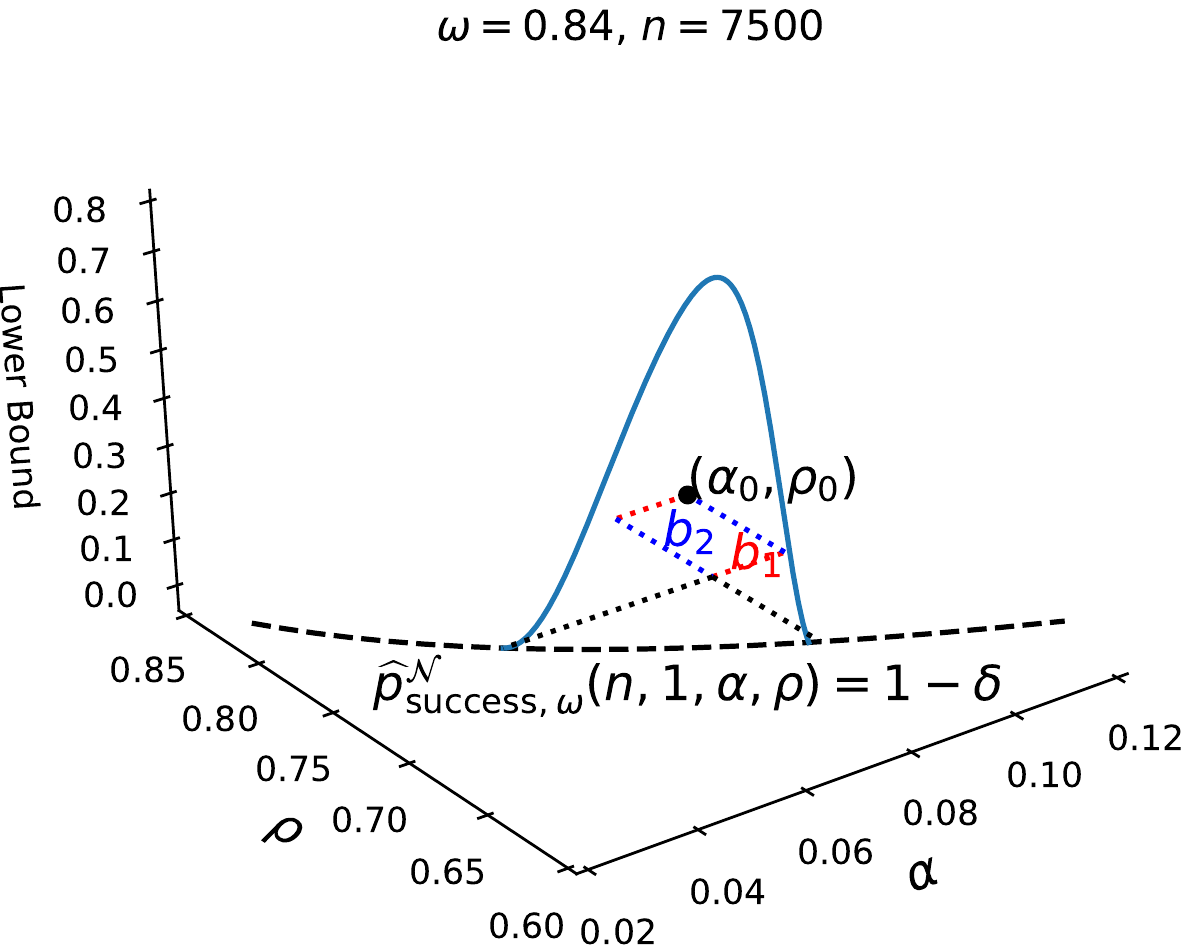}\centering
\caption{A visual depiction of how the bound \eqref{eq:lb_stopping_optimised} is optimised. The solid curve is a plot of the original bound \eqref{eq:lb_stopping} along the Pareto front (dashed curve), for given $\omega = 0.84$, and with $n = 7500$, $\alpha_{0} = 0.1$, $\rho_{0} = 0.8$, $\delta = 0.1$, $\beta_{1} = 0.05$, $\beta_{2} = 0.05$. The dotted lines indicate how the gaps $\left(\alpha_{0}-b_{1}\right)-\alpha^{*}$ and $\left(\rho_{0}-b_{2}\right) -\rho^{*}$ must be positive in order for the bound to be informative.}
\label{fig:optimised}
\end{figure}

\bibliography{autosam_web}{} 

\begin{thebibliography}{10}

\bibitem{Tempo2013}
R.~Tempo, G.~Calafiore, and F.~Dabbene, {\em Randomized Algorithms for Analysis
  and Control of Uncertain Systems With Applications}.
\newblock Springer, 2nd~ed., 2013.

\bibitem{Khargonekar1996}
P.~Khargonekar and A.~Tikku, ``Randomized algorithms for robust control
  analysis and synthesis have polynomial complexity,'' in {\em 35th {IEEE}
  Conference on Decision and Control}, {IEEE}, 1996.

\bibitem{Tempo1996}
R.~Tempo, E.~W. Bai, and F.~Dabbene, ``Probabilistic robustness analysis:
  explicit bounds for the minimum number of samples,'' in {\em 35th {IEEE}
  Conference on Decision and Control}, {IEEE}, 1996.

\bibitem{Vidyasagar1998}
M.~Vidyasagar, ``Statistical learning theory and randomized algorithms for
  control,'' {\em {IEEE} Control Systems Magazine}, vol.~18, no.~6, pp.~69--85,
  1998.

\bibitem{Ding2019}
S.~X. Ding, L.~Li, and M.~Kr{\"{u}}ger, ``Application of randomized algorithms
  to assessment and design of observer-based fault detection systems,'' {\em
  Automatica}, vol.~107, pp.~175--182, 2019.

\bibitem{Alpcan2003}
T.~Alpcan, T.~Basar, and R.~Tempo, ``Randomized algorithms for stability and
  robustness analysis of high speed communication networks,'' in {\em {IEEE}
  Conference on Control Applications}, {IEEE}, 2003.

\bibitem{Ho2007}
Y.-C. Ho, Q.-C. Zhao, and Q.-S. Jia, {\em Ordinal Optimization: Soft
  Optimization for Hard Problems}.
\newblock Springer, 2007.

\bibitem{Ho1992}
Y.~C. Ho, R.~S. Sreenivas, and P.~Vakili, ``Ordinal optimization of {DEDS},''
  {\em Discrete Event Dynamic Systems}, vol.~2, no.~1, pp.~61--88, 1992.

\bibitem{Xie1997}
X.~Xie, ``Dynamics and convergence rate of ordinal comparison of stochastic
  discrete-event systems,'' {\em {IEEE} Transactions on Automatic Control},
  vol.~42, no.~4, pp.~586--590, 1997.

\bibitem{Lee1999}
L.~H. Lee, T.~W.~E. Lau, and Y.~C. Ho, ``Explanation of goal softening in
  ordinal optimization,'' {\em {IEEE} Transactions on Automatic Control},
  vol.~44, no.~1, pp.~94--99, 1999.

\bibitem{Deng1999}
M.~Deng and Y.-C. Ho, ``An ordinal optimization approach to optimal control
  problems,'' {\em Automatica}, vol.~35, no.~2, pp.~331--338, 1999.

\bibitem{Ho1995}
Y.-C. Ho and M.~E. Larson, ``Ordinal optimization approach to rare event
  probability problems,'' {\em Discrete Event Dynamic Systems: Theory and
  Applications}, vol.~5, no.~2-3, pp.~281--301, 1995.

\bibitem{Yang2002}
M.~Yang and L.~Lee, ``Ordinal optimization with subset selection rule,'' {\em
  Journal of Optimization Theory and Applications}, vol.~113, no.~3,
  pp.~597--620, 2002.

\bibitem{Vidyasagar2001}
M.~Vidyasagar, ``Randomized algorithms for robust controller synthesis using
  statistical learning theory,'' {\em Automatica}, vol.~37, no.~10,
  pp.~1515--1528, 2001.

\bibitem{Ishii2008}
H.~Ishii and R.~Tempo, ``Las vegas randomized algorithms in distributed
  consensus problems,'' in {\em American Control Conference}, {IEEE}, 2008.

\bibitem{Chin2021}
R.~Chin, J.~E. Rowe, I.~Shames, C.~Manzie, and D.~Ne{\v{s}}i{\'{c}}, ``Ordinal
  optimisation and the offline multiple noisy secretary problem.''
  https://arxiv.org/abs/2106.01185, 2021.

\bibitem{Joe2014}
H.~Joe, {\em Dependence Modeling with Copulas}.
\newblock CRC Press, 2014.

\bibitem{Dubhashi2009}
D.~P. Dubhashi and A.~Panconesi, {\em Concentration of Measure for the Analysis
  of Randomized Algorithms}.
\newblock Cambridge University Press, 2009.

\bibitem{Calafiore2006}
G.~Calafiore and M.~Campi, ``The scenario approach to robust control design,''
  {\em {IEEE} Transactions on Automatic Control}, vol.~51, no.~5, pp.~742--753,
  2006.

\bibitem{Koltchinskii2000}
V.~Koltchinskii, C.~T. Abdallah, M.~Ariola, P.~Dorato, and D.~Panchenko,
  ``Improved sample complexity estimates for statistical learning control of
  uncertain systems,'' {\em {IEEE} Transactions on Automatic Control}, vol.~45,
  no.~12, pp.~2383--2388, 2000.

\bibitem{Fujisaki2007}
Y.~Fujisaki and Y.~Oishi, ``Guaranteed cost regulator design: A probabilistic
  solution and a randomized algorithm,'' {\em Automatica}, vol.~43, no.~2,
  pp.~317--324, 2007.

\bibitem{Alamo2015}
T.~Alamo, R.~Tempo, A.~Luque, and D.~R. Ramirez, ``Randomized methods for
  design of uncertain systems: Sample complexity and sequential algorithms,''
  {\em Automatica}, vol.~52, pp.~160--172, 2015.

\bibitem{Bayraksan2012}
G.~Bayraksan and P.~Pierre-Louis, ``Fixed-width sequential stopping rules for a
  class of stochastic programs,'' {\em {SIAM} Journal on Optimization},
  vol.~22, no.~4, pp.~1518--1548, 2012.

\bibitem{Grammatico2015}
S.~Grammatico, X.~Zhang, K.~Margellos, P.~Goulart, and J.~Lygeros, ``A scenario
  approach for non-convex control design,'' {\em {IEEE} Transactions on
  Automatic Control}, pp.~1--1, 2015.

\bibitem{Esfahani2015}
P.~M. Esfahani, T.~Sutter, and J.~Lygeros, ``Performance bounds for the
  scenario approach and an extension to a class of non-convex programs,'' {\em
  {IEEE} Transactions on Automatic Control}, vol.~60, no.~1, pp.~46--58, 2015.

\bibitem{Luedtke2008}
J.~Luedtke and S.~Ahmed, ``A sample approximation approach for optimization
  with probabilistic constraints,'' {\em {SIAM} Journal on Optimization},
  vol.~19, no.~2, pp.~674--699, 2008.

\bibitem{Kendall1990}
M.~Kendall and J.~D. Gibbons, {\em Rank Correlation Methods}.
\newblock Oxford University Press, 5th~ed., 1990.

\bibitem{Gray2011}
R.~M. Gray, {\em Entropy and Information Theory}.
\newblock Springer, 2011.

\bibitem{Liu2012a}
H.~Liu, F.~Han, M.~Yuan, J.~Lafferty, and L.~Wasserman, ``High-dimensional
  semiparametric gaussian copula graphical models,'' {\em The Annals of
  Statistics}, vol.~40, no.~4, pp.~2293--2326, 2012.

\bibitem{Arnold2008}
B.~C. Arnold, N.~Balakrishnan, and H.~N. Nagaraja, {\em A First Course in Order
  Statistics}.
\newblock SIAM, 2008.

\bibitem{Capinski2004}
M.~Capinski and P.~E. Kopp, {\em Measure, Integral and Probability}.
\newblock Springer, 2004.

\bibitem{Ono2008}
M.~Ono and B.~C. Williams, ``Iterative risk allocation: A new approach to
  robust model predictive control with a joint chance constraint,'' in {\em
  47th {IEEE} Conference on Decision and Control}, {IEEE}, 2008.

\bibitem{Shekhar2017}
R.~C. Shekhar, G.~S. Sankar, C.~Manzie, and H.~Nakada, ``Efficient calibration
  of real-time model-based controllers for diesel engines {\textemdash} part i:
  Approach and drive cycle results,'' in {\em {IEEE} 56th Annual Conference on
  Decision and Control ({CDC})}, {IEEE}, 2017.

\bibitem{Maass2020}
A.~I. Maass, C.~Manzie, I.~Shames, R.~Chin, D.~Ne\v{s}i\'{c}, N.~Ulapane, and
  H.~Nakada, ``Tuning of model predictive engine controllers over transient
  drive cycles,'' in {\em 21st IFAC World Congress}, 2020.

\bibitem{Sankar2020}
G.~S. Sankar, R.~C. Shekhar, C.~Manzie, T.~Sano, and H.~Nakada, ``Fast
  calibration of a robust model predictive controller for diesel engine
  airpath,'' {\em {IEEE} Transactions on Control Systems Technology}, vol.~28,
  pp.~1505--1519, jul 2020.

\bibitem{Chin2020a}
R.~Chin, A.~I. Maass, N.~Ulapane, C.~Manzie, I.~Shames, D.~Ne{\v{s}}i{\'{c}},
  J.~E. Rowe, and H.~Nakada, ``Active learning for linear parameter-varying
  system identification,'' in {\em 21th IFAC World Congress}, 2020.

\bibitem{Ira2020}
A.~S. Ira, C.~Manzie, I.~Shames, R.~Chin, D.~Ne{\v{s}}i{\'{c}}, H.~Nakada, and
  T.~Sano, ``Tuning of multivariable model predictive controllers through
  expert bandit feedback,'' {\em International Journal of Control}, 2020.

\bibitem{Bear}
``Birmingham environment for academic research ({BEAR}).''
  http://www.birmingham.ac.uk/bear.

\bibitem{Lafayette2016}
L.~Lafayette, G.~Sauter, L.~Vu, and B.~Meade, ``Spartan performance and
  flexibility: An hpc-cloud chimera,'' in {\em OpenStack Summit}, 2016.

\bibitem{Hoeffding1963}
W.~Hoeffding, ``Probability inequalities for sums of bounded random
  variables,'' {\em Journal of the American Statistical Association}, vol.~58,
  no.~301, pp.~13--30, 1963.

\end{thebibliography}
\bibliographystyle{ieeetr}

\end{document}